\documentclass[12pt]{amsart}
\usepackage{amssymb,amsmath,eucal}
\usepackage[dvips]{graphics}

\newcommand{\la}{\lambda}
\newcommand{\al}{\alpha}
\newcommand{\be}{\beta}
\newcommand{\Ga}{\Gamma}

\newcommand{\Prob}{\operatorname{Prob}}
\newcommand{\Conf}{\operatorname{Conf}}
\newcommand{\CX}{\Conf(\X)}
\newcommand{\CXf}{\CX_0}
\newcommand{\sgn}{\operatorname{sgn}}
\newcommand{\tr}{\operatorname{tr}}
\newcommand{\tri}{\triangle}
\newcommand{\trZ}{\tri_Z}
\newcommand{\X}{\mathfrak{X}}

\newcommand{\A}{\mathsf{A}}
\newcommand{\J}{\mathsf{J}}

\newcommand{\wt}{\widetilde}
\newcommand{\sD}{\mathcal{D}}
\newcommand{\brh}{\boldsymbol{\varrho}}

\newcommand{\Z}{\mathbb{Z}}
\newcommand{\R}{\mathbb{R}}
\newcommand{\C}{\mathbb{C}}
\newcommand{\Fr}{\operatorname{Fr}}
\newcommand{\const}{\operatorname{const}}
\newcommand{\bS}{\mathsf{S}}
\newcommand{\h}{\frac12}
\newcommand{\sh}{{\textstyle \frac12}}
\newcommand{\KK}{\mathsf{K}}

\newcommand{\LL}{\mathsf{L}}
\newcommand{\DD}{\mathsf{D}}
\newcommand{\kk}{\mathsf{k}} 
\newcommand{\De}{\Delta} 
\newcommand{\bo}{\bar 1}

\newcommand{\ep}{\varepsilon}
\newcommand{\Lim}{\operatorname{\mathcal{ L}im}}
\newcommand{\g}{\gamma}
\newcommand{\aF}{\mathcal{F}}

\newcommand{\ix}{\theta}
\newcommand{\Lc}{\mathcal{L}}
\newcommand{\LH}{\mathcal{L}_{1|2}(H)}

\theoremstyle{plain}
\newtheorem{theorem}{Theorem}
\newtheorem{lemma}{Lemma}[section]
\newtheorem{proposition}[lemma]{Proposition}
\newtheorem{corollary}[lemma]{Corollary}

\theoremstyle{definition}
\newtheorem{remark}[lemma]{Remark}
\newtheorem{definition}[lemma]{Definition}

\numberwithin{equation}{section}

\begin{document}

\title[Asymptotics of Plancherel Measures]
{Asymptotics of Plancherel measures for symmetric groups} 
\author[A.~Borodin]{Alexei Borodin}
\address{University of Pennsylvania, Department of Mathematics,
Philadelphia PA 19104--6395 }
\email{borodine@math.upenn.edu}
\author[A.~Okounkov]{Andrei Okounkov}
\address{University of Chicago, Department of Mathematics,
5734 University Ave., Chicago IL 60637}
\email{okounkov@math.berkeley.edu}
\curraddr{Department of Mathematics, University of California
at Berkeley, Evans Hall, Berkeley, CA 94720-3840}
\author[G.~Olshanski]{Grigori Olshanski}
\address{Dobrushin Mathematics Laboratory,
Institute for Problems of Information Transmission, Bolshoy Karetny 19,
101447, Moscow, Russia}
\thanks{A.O.\ is supported by NSF grant DMS-9801466, G.O.\ is supported by the Russian      
Foundation for Basic Research under grant 98-01-00303}
\email{olsh@glasnet.ru}
\begin{abstract}
We consider the asymptotics of the Plancherel measures
on partitions of $n$ as $n$ goes to infinity. We prove
that the local structure of a Plancherel typical partition
in the middle of the limit shape converges to a determinantal
point process with the discrete sine kernel. 

On the edges
of the limit shape, we prove that the joint distribution
of suitably scaled 1st, 2nd, and so on rows of a Plancherel
typical diagram converges to the corresponding distribution for
eigenvalues of random Hermitian matrices (given by the Airy
kernel). This
proves a conjecture due to Baik, Deift, and Johansson by
methods different from the Riemann-Hilbert techniques used in their
original
papers \cite{BDJ1,BDJ2} and from the combinatorial
proof given in \cite{O}.
 
Our approach is based on an exact determinantal
formula  for the correlation functions of the
poissonized Plancherel measures in terms of a new 
kernel involving Bessel functions.  Our
asymptotic analysis relies on the
classical asymptotic formulas for the Bessel functions and
depoissonization techniques.
\end{abstract}

\maketitle 

\section{Introduction}\label{s1}

\subsection{Plancherel measures}\label{s11}
Given  a finite group $G$,  by the corresponding Plancherel
measure we mean 
the probability  measure on the set $G^\wedge$ of irreducible
representations of $G$ which assigns to a representation 
$\pi\in G^\wedge$ the weight $(\dim\pi)^2/|G|$.
For the symmetric group $S(n)$, the set $S(n)^\wedge$ is
the set of partitions $\la$ of the number $n$, which
we shall identify with Young diagrams with $n$ squares
throughout this paper.  
The Plancherel measure on partitions $\la$
 arises naturally in representation--theoretic,
combinatorial,  and probabilistic problems. For
example, the Plancherel distribution of the first part
of a partition coincides with the distribution
of the longest increasing subsequence of a
uniformly distributed random permutation \cite{Sch}. 

We denote the Plancherel measure
on partitions of $n$ by $M_n$
\begin{equation}\label{e0001}
M_n(\la)=\frac{(\dim\la)^2}{n!} \,, \quad |\la|=n \,,
\end{equation}
where $\dim\la$ is the dimension of the corresponding
representation of $S(n)$. The asymptotic properties of these
measures as $n\to\infty$ have been studied very intensively, see the 
References and below. 

\begin{figure}[!h]
\centering
\scalebox{.85}{\includegraphics{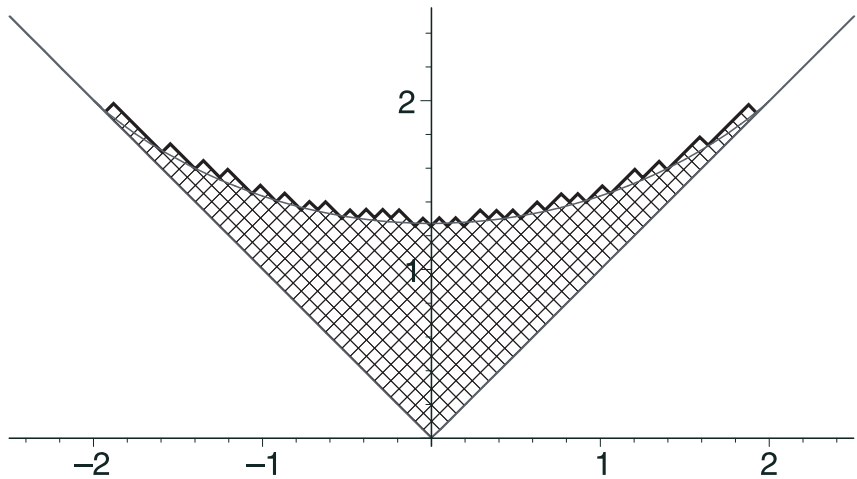}}
\caption{The limit shape of a typical diagram.}
\label{fig1}
\end{figure}

In the seventies, Logan and Shepp \cite{LS} and,
independently, Vershik and Kerov \cite{VK1,VK2} discovered
the following measure concentration phenomenon for
$M_n$ as $n\to\infty$.   Let $\la$ be a partition of $n$ and let $i$ and $j$ be the usual 
coordinates on the diagrams, namely, the row number and the column number.
Introduce new coordinates $u$ and $v$ by 
\begin{equation*}
u=\frac{j-i}{\sqrt{n}} \,, \quad  v=\frac{i+j}{\sqrt{n}} \,,
\end{equation*}
that is, we flip the diagram, rotate it $135^\circ$ as in Figure \ref{fig1},
and scale it by the factor of $n^{-1/2}$ in both directions.

After this scaling, the Plancherel measures $M_n$    
converge as $n\to\infty$ (see \cite{LS,VK1,VK2} for precise statements) to the delta measure 
supported on the following shape:
\begin{equation*}
\{ |u|\le 2, |u| \le v \le \Omega(u)\}\,,
\end{equation*}
where the function $\Omega(u)$ is defined by
\begin{equation*}
\Omega(u)=\cases
{\frac2\pi\, \left(u \arcsin(u/2) + \sqrt{4-u^2}\right)}\,, & |u|\le 2 \,,\\
|u|\,, & |u|>2 \,.
\endcases 
\end{equation*}
The function $\Omega$ is plotted in Figure \ref{fig1}. 
As explained in great detail in \cite{Ke5}, this limit shape $\Omega$ is
very closely connected to Wigner's semicircle law for distribution of
eigenvalues of a random matrices, see also \cite{Ke2,Ke3,Ke4}. 

{}From a different point of view, the connection with random
matrices was observed in \cite{BDJ1,BDJ2}, and also in the 
earlier papers \cite{J,R,Re}. In \cite{BDJ1}, Baik, Deift, and
Johansson made the following conjecture. They conjectured that
in the $n\to\infty$ limit and after proper scaling the 
joint distribution of $\la_i$, $i=1,2,\dots$ becomes
identical to the joint distribution of largest eigenvalues
of a Gaussian random Hermitian matrix (which is known to be
the so-called Airy ensemble, see Section \ref{s14}). They 
proved this for the  individual distribution
of $\la_1$ and $\la_2$ in \cite{BDJ1} and \cite{BDJ2}, respectively. 
A combinatorial proof of the full conjecture was given by one
of us in \cite{O}. It was based on an interplay
between maps on surfaces and ramified coverings
of the sphere. 

In this paper we study the local structure of a typical Plancherel
diagram both in the bulk of the limit shape $\Omega$ and on its edge,
where by the study of the edge we mean the study of the behavior
of $\la_1$, $\la_2$, and so on.

We employ an analytic approach 
based on an exact formula in terms of Bessel functions
for the correlation functions of the
so-called \emph{poissonization} of the Plancherel measures $M_n$,
see Theorem \ref{t1} in the  following section, and the so-called \emph{depoissonization}
techniques, see Section 1.4. 

The exact formula in Theorem \ref{t1}
is a limit case of a formula from \cite{BO}, see
also the recent paper \cite{O2} for a more general result. 
The use of poissonization and 
depoissonization is very much in the spirit of \cite{BDJ1,J,V}
and represents a well known in statistical mechanics principle of
the equivalence of ensembles.

Our main results are the following two. In the bulk of the limit
shape $\Omega$, we prove that the local
structure of a Plancherel typical partition converges to  
a determinantal  point process with the discrete sine kernel,
see Theorem \ref{t2}. This result is parallel to the corresponding result for
random matrices.  On the edge of the limit shape, we give
an analytic proof of  the Baik-Deift-Johansson conjecture,
see Theorem \ref{t3}. These results will be stated in
Subsections \ref{s13} and \ref{s14} of the present
Introduction, respectively. 

Simultaneously and independently, results equivalent to our
Theorems \ref{t1b} and \ref{t3} were obtained by K.~Johansson
\cite{J2}. 

\subsection{Poissonization and correlation functions}\label{s12}

For $\ix>0$, consider 
the {\it poissonization} $M^\ix$ of the measures $M_n$
\begin{equation*}
M^\ix (\la) = e^{-\ix}\sum_n \frac{\ix^n}{n!}\, M_n(\la) =
e^{-\ix} \ix^{|\la|} \,
\left(\frac{\dim\la}{|\la|!}
\right)^2 \,. 
\end{equation*}
This is a probability measure on the set of all partitions. 
Our first result is the computation of the correlation
functions of the measures $M^\ix$. 

By correlation functions we mean the following. 
By definition, set
$$
\sD(\la)=\{\la_i-i\}\subset\Z \,.
$$
Also, following \cite{VK3}, define  the {\em modified 
Frobenius coordinates} $\Fr(\la)$ of a partition $\la$ by 
\begin{multline}\label{e0009}
%\la\mapsto
\Fr(\la)=\left(\sD(\la)+\sh\right)\triangle \left(\Z_{\le 0} -\sh\right)\\
=
\left\{p_1+\sh,\dots,p_d+\sh,-q_1-\sh,\dots,-q_d-\sh
\right\} \subset \Z+\sh\,,
\end{multline}
where $\triangle$ stands for the symmetric difference of two sets,
$d$ is the number of squares on the diagonal of $\la$, 
and $p_i$'s and $q_i$'s are the usual Frobenius
coordinates of $\la$. Recall that $p_i$ is the number of
squares in the $i$th row  to the right  of the diagonal,
and $q_i$ is number of squares in the $i$th column
below the diagonal. The equality \eqref{e0009} is a well known
combinatorial fact discovered by
Frobenius, see Ex.~I.1.15(a) in \cite{M}.
Note that, in contrast to $\Fr(\la)$, the
set $\sD(\la)$ is infinite and, moreover, it contains all but
finitely many negative integers. 

The sets $\sD(\la)$ and $\Fr(\la)$ have the following nice
geometric interpretation. Let the diagram $\la$ be flipped and
rotated $135^\circ$ as in Figure \ref{fig1}, but not scaled.
Denote by $\omega_\la$
a piecewise linear function with $\omega'_\la=\pm 1$ whose
graph if given by the upper boundary of $\la$ completed by the lines
$$
v=|u|\,, \quad u\notin[-\la'_1,\la_1] \,.
$$
Then 
$$
k\in\sD(\la) \Leftrightarrow \omega'_\la\Big|_{[k,k+1]}=-1 \,.
$$
In other words, if we consider $\omega_\la$ as a history
of a walk on $\Z$ then $\sD(\la)$ are those moments when
a step is made in the negative direction. It is
therefore natural to call $\sD(\la)$ the \emph{descent set}
of $\la$. As we shall see, the correspondence $\la\mapsto\sD(\la)$
is a very convenient way to encode the local structure of the boundary of $\la$. 

The halves in the definition of 
$\Fr(\la)$ have the following  interpretation:
one  splits the diagonal squares in half and gives half to the rows
and half to the columns.  

\begin{definition} The correlation functions of $M^\ix$ are the probabilities
that the sets $\Fr(\la)$ or, similarly, $\sD(\la)$ contain a fixed
subset $X$. More precisely, we set
\begin{alignat}{2}
\rho^\ix(X)&= M^\ix\left(\left\{\la\,|\, X\subset \Fr(\la)\right\} 
\right)\,,\quad &&X\in\Z+\sh\,,\label{e11a}\\
\brh^\ix(X)&= M^\ix\left(\left\{\la\,|\, X\subset \sD(\la)\right\} 
\right)\,,\quad &&X\in\Z\,. \label{e0006}
\end{alignat}
\end{definition}

\begin{theorem}\label{t1}
For any $X=\{x_1,\dots,x_s\}\subset \Z+\sh$
we have 
\begin{equation*}
\rho^\ix(X)=
\det \Big[\KK(x_i,x_j)\Big]_{1\le i,j \le s}\,,
\end{equation*}
where the kernel $\KK$ is given by the following formula
\begin{equation}\label{e12a}
\KK(x,y)=\left\{
\begin{array}{ll}
\displaystyle
\sqrt{\ix}\,\,\frac{\kk_+(|x|,|y|)}{|x|-|y|}\,, & xy >0 \,,\\[10pt]
\displaystyle \sqrt{\ix}\,\,\frac{\kk_-(|x|,|y|)}{x-y} \,, & xy<0 \,.
\end{array}
\right. 
\end{equation}
The functions $\kk_\pm$ are defined by
\begin{align}\label{e12} 
\kk_+(x,y)&= 
J_{x-\h}\, J_{y+\h}  - J_{x+\h} \, J_{y-\h} \,, \\ 
\kk_-(x,y)&= J_{x-\h}\, J_{y-\h}  + J_{x+\h} \, J_{y+\h} \label{e13}\,,
\end{align}
where $J_x=J_x(2\sqrt{\ix})$ is the Bessel function of
order $x$ and argument $2\sqrt{\ix}$. 
\end{theorem}

This theorem is established in Section \ref{s21}, see also
Remark \ref{r11} below.  By the complementation principle,
see Sections \ref{sA3} and \ref{s21b}, Theorem \ref{t1} is
equivalent to the following 

\begin{theorem}\label{t1b}
For any $X=\{x_1,\dots,x_s\}\subset \Z$
we have 
\begin{equation}\label{e0005}
\brh^\ix(X)=
\det \Big[\J(x_i,x_j)\Big]_{1\le i,j \le s}\,,
\end{equation}
Here the kernel $\J$ is given by the following formula
\begin{equation}\label{e212}
\J(x,y)=\J(x,y;\ix)=\sqrt{\ix}\,\, \frac{J_x\, J_{y+1} - J_{x+1} \, J_{y}}
{x-y} \,,
\end{equation}
where $J_x=J_x(2\sqrt\ix)$. 
\end{theorem}

\begin{remark}\label{r11}
Theorem \ref{t1} is a limit case of Theorem 3.3  of \cite{BO}.
For the reader's convenience a direct proof of it is given
in Section \ref{s2}.  Another proof of the results of \cite{BO}
will appear in \cite{BO2}. Various limit cases of the results
of \cite{BO} are discussed in \cite{BO3}. 
By different methods, the formula \eqref{e0005} was obtained by
K.~Johansson \cite{J2}.

A representation--theoretic
proof of a more general formula than Theorem 3.3 of \cite{BO}
has been subsequently  given in \cite{O2}. 
\end{remark}

\begin{remark}
Observe that all Bessel functions involved 
in the above formulas  are of integer order.
Also note that the ratios like  $\J(x,y)$ are entire functions of $x$ and $y$
because  $J_x$ is an entire function of $x$. In particular, the values $\J(x,x)$
are well defined. Various
denominator--free formulas for the kernel $\J$ are given  in Section \ref{s21}.
\end{remark}

\subsection{Asymptotics in the bulk of the spectrum}\label{s13} 

Given a sequence of subsets
\begin{equation*}
X(n)=\left\{x_1(n)<\dots<x_s(n)\right\}\subset\Z\,, 
\end{equation*}
where $s=|X(n)|$ is  some fixed integer, we call this sequence {\em regular} 
 if the following limits 
\begin{align}
a_i &= \lim_{n\to\infty} \frac{x_i(n)}{\sqrt{n}}\,, \label{e14}\\
d_{ij}&=\lim_{n\to\infty}\left(x_i(n)-x_j(n)\right)\,, \label{e15}
\end{align}
exist, finite or infinite. Here $i,j=1,\dots,s$. Observe that if $d_{ij}$ is finite then
$d_{ij}=x_i(n)-x_j(n)$ for $n\gg 0$.

In the case when $X(n)$ can be represented
as $X(n)=X'(n)\cup X''(n)$ and the distance between $X'(n)$ and
$X''(n)$ goes to $\infty$ as $n\to\infty$ we shall say that the
sequence {\it  splits}; otherwise, we call it {\em nonsplit}.
Obviously, $X(n)$ is nonsplit if and only if all $x_i(n)$ stay at a finite
distance from each other.

Define the correlation
functions $\brh(n,\,\cdot\,)$ of the measures $M_n$ by
the same rule as in \eqref{e0006}
\begin{equation*}
\brh(n,X)= M_n\left(\left\{\la\,|\, X\subset \sD(\la)\right\} 
\right)\,.
\end{equation*}
We are interested in the limit of $\brh(n,X(n))$ as $n\to\infty$.
This limit will be computed in
Theorem \ref{t2} below. As we shall see, if $X(n)$ splits, then
the limit correlations factor accordingly. 

Introduce the following {\em discrete sine kernel}  which is
a translation invariant kernel on the lattice $\Z$
$$
\bS(k,l;a)=\bS(k-l,a)\,, \quad k,l\in\Z\,,
$$
depending on a real parameter $a$: 
\begin{align*}
\bS(k,a)&=\frac{\sin(\arccos(a/2)\,k)}{\pi k}\\
&=
\frac{\sqrt{4-a^2}}{2\pi} \frac{U_{k-1}(a/2)}{k}\,,
\quad k\in \Z \,.
\end{align*}
Here $U_k$ is 
the Tchebyshev polynomials of the second kind.
We agree that 
\begin{equation*}
\bS(0,a)=\frac{\arccos(a/2)}\pi\,, 
\quad \bS(\infty,a)=0
\end{equation*}
and also that 
$$
\bS(k,a)=
\cases 0\,,& \textup{$a\ge 2$ or $a\le 2$ and $k\ne 0$}\,,\\
1 \,,  &\textup{$a\le 2$ and $k= 0$}\,.
\endcases
$$

The following result describes the local structure
of a Plancherel typical partition. 

\begin{theorem}\label{t2}
 Let $X(n)\subset\Z$ be a regular sequence
and let the numbers $a_i$, $d_{ij}$ be defined by \eqref{e14}, \eqref{e15}.  
If $X(n)$ splits, that is, if $X(n)=X'(n)\cup X''(n)$ and the distance between $X'(n)$ and
$X''(n)$ goes to $\infty$ as $n\to\infty$ then 
\begin{equation}\label{e1t2}
\lim_{n\to\infty} \brh(n,X(n)) =
\lim_{n\to\infty} \brh(n,X'(n))  \cdot
\lim_{n\to\infty} \brh(n,X''(n))\,.  
\end{equation}
If $X(n)$ is nonsplit  then 
\begin{equation}\label{e2t2}
\lim_{n\to\infty} \brh\left(n,X(n)\right)
=
\det \Big[\,\bS(d_{ij},a)
\Big]_{1\le i,j \le s}\,,
\end{equation}
where $\bS$ is the discrete sine kernel and $a=a_1=a_2=\dots$. 
\end{theorem}

We prove this theorem in Section \ref{s3}. 

\begin{remark} Notice that, in particular, Theorem \ref{t2}
implies that, as $n\to\infty$, the shape of a typical partition
$\la$ near any point of the limit curve $\Omega$ is described
by a stationary random process. For distinct points on the
curve $\Omega$ these random processes are independent. 
\end{remark}

\begin{remark} By complementation, see
Section \ref{sA3} and \ref{s32}, one obtains from Theorem \ref{t2}
an equivalent statement about the 
asymptotics of the following correlation functions
\begin{equation*}
\rho(n,X)= M_n\left(\left\{\la\,|\, X\subset \Fr(\la)\right\} 
\right)\,.
\end{equation*}
\end{remark} 

\begin{remark} The discrete sine kernel was studied before, see \cite{Wi1,Wi2},
mainly as a model case for the continuous sine kernel. In particular, the
asymptotics of Toeplitz determinants built from the discrete sign kernel
was obtained by H.~Widom in \cite{Wi2} answering a question of F.~Dyson. As pointed
out by S.~Kerov, this asymptotics has interesting consequences for the Plancherel
measures. 
\end{remark} 

\begin{remark}
Note that, in particular, Theorem \ref{t2} implies that the limit
density (the 1-point
correlation function)  is given by  
\begin{equation}\label{e16}
\brh(\infty,a)=\cases
\frac1\pi\, {\arccos(a/2)}\,, & |a|\le 2\,,\\
0\,, & a>2 \,,\\
1\,, & a<-2 \,.
\endcases 
\end{equation}
This is in agreement with the Logan-Shepp-Vershik-Kerov result
about the limit shape $\Omega$. More concretely,
the function $\Omega$ is related to the density \eqref{e16} by
\begin{equation*}
\brh(\infty,u)=\frac{1-\Omega'(u)}2 \,,
\end{equation*}
which can be interpreted as follows. Approximately, 
we have 
\begin{equation*}
\#\left\{i\,\left|\, \frac{\la_i}{\sqrt{n}} \in [u,u+\De u]\right.\right\}
\approx \sqrt{n} \, \brh(\infty,u) \, \De u \,.
\end{equation*}
Set $w=\dfrac{i}{\sqrt{n}}$. Then the above relation reads 
$\De w \approx \brh(\infty,u) \, \De u $ and it should be 
satisfied on the boundary  $v=\Omega(u)$ of the limit shape.
Since $v=u+2w$, we conclude that
\begin{equation*}
\brh(\infty,u)\approx \frac{d w}{du}=\frac{1-\Omega'}2 \,,
\end{equation*}
as was to be shown.
\end{remark}

\begin{remark}
The discrete sine-kernel $\bS$ becomes especially
nice near the diagonal, that is, where $a=0$. Indeed,
\begin{equation*}
\bS(x,0)= \cases 
1/2\,, & x=0 \,,\\
{(-1)^{(x-1)/2}}\Big/(\pi x) \,,
& x=\pm1,\pm3,\dots \,, \\
0\,, &x=\pm2,\pm4,\dots \,.
\endcases
\end{equation*}
\end{remark}

\subsection{Behavior near the edge of the spectrum and the Airy ensemble}\label{s14}

The discrete sine kernel $\bS(k,a)$ vanishes if $a\ge 2$. Therefore,
it follows from Theorem \ref{t2} that the limit correlations $\lim\brh(n,X(n))$
vanish if $a_i\ge 2$ for some $i$.  However, as will be shown below in
Proposition \ref{p41}, 
after a suitable scaling near the edge $u=2$,  the
correlation functions $\brh^\ix$ converge to 
the correlation functions given
by the Airy kernel \cite{F,TW}
\begin{equation*}
\A(x,y)=\frac{A(x)A'(y)-A'(x)A(y)}{x-y}\,.  
\end{equation*}
Here $A(x)$ is the Airy function: 
\begin{equation}\label{e17}
A(x)=\frac1\pi\,\int_0^\infty \cos\left(\frac{u^3}3+xu\right)\,du.  
\end{equation}
In fact,  the following  more precise 
statement is true about the behavior
of the Plancherel measure near the edge $u=2$. By
symmetry, everything we say about the edge 
$u=2$ applies to the opposite edge $u=-2$. 

Consider the random point process on $\R$ whose correlation
functions are given by the determinants 
\begin{equation*}
\rho^{\textup{Airy}}_k(x_1,\dots,x_k)= \det \Big[\,\A(x_i,x_j)
\Big]_{1\le i,j \le k}
\end{equation*}
and let 
\begin{equation*}
\zeta=(\zeta_1 > \zeta_2 > \zeta_3 > \dots ) \in\R^\infty
\end{equation*}
be its random configuration. We call the random
variables $\zeta_i$'s the {\it Airy ensemble}. 
It is known \cite{F,TW} that the Airy ensemble
describes the 
behavior of the (properly scaled) 1st, 2nd, and so
on largest eigenvalues of a Gaussian random
Hermitian matrix. The distribution
of individual eigenvalues was obtained by
Tracy and Widom in \cite{TW} in terms of
certain Painlev\'e transcendents.

It has been conjectured by
Baik, Deift, and Johansson that 
the random variables
\begin{equation*}
\widetilde{\la}=
\left(\widetilde{\la}_1 \ge \widetilde{\la}_2 \ge \dots
\right)\,,
\quad 
\widetilde{\la}_i =
n^{1/3} \, \left(
\frac{\la_i}{n^{1/2}}-2 \right) 
\end{equation*}
converge, in distribution and together with all moments,
to the Airy ensemble. They verified this
conjecture for individual distribution of $\la_1$ and 
$\la_2$ in \cite{BDJ1} and \cite{BDJ2}, respectively. 
In particular, in the case of $\la_1$, this generalizes the
result of  \cite{VK1,VK2} that $\frac{\la_1}{\sqrt n} \to 2$
in probability as $n\to\infty$. The computation of $\lim \frac{\la_1}{\sqrt n}$ was
known as the Ulam problem;  different  solutions to this problem were
given in \cite{AD,J,S}.

Convergence of all expectations
of the form
\begin{equation}\label{e18a}
\left\langle
\prod_{k=1}^r 
\sum_{i=1}^\infty e^{t_k \widetilde{\la}_i}
\right\rangle\,,
\quad t_1,\dots,t_r>0\,, 
\quad r=1,2,\dots\,,
\end{equation}
to the corresponding quantities 
for the Airy ensembles was established in \cite{O}. The proof 
in \cite{O} was based on a combinatorial interpretation of
\eqref{e18a} as the asymptotics in a certain enumeration problem
for random surfaces. 

In the present paper we use different ideas to prove the following 

\begin{theorem} \label{t3}  As $n\to\infty$, the random 
variables $\widetilde{\la}$ 
converge, in joint distribution, to the Airy ensemble.
\end{theorem}

This is done in Section \ref{s4} using methods described in the next
subsection. The result stated in Theorem \ref{t3} was
independently obtained by K.~Johansson in \cite{J2}. 

\subsection{Poissonization and depoissonization}\label{s15}

We obtain Theorems \ref{t2} and \ref{t3} from Theorem \ref{t1} using the
so-called depoissonization techniques. We recall that  the fundamental idea of 
depoissonization is the following. 

Given a sequence $b_1,b_2,b_3,\dots$ its {\em poissonization}
is, by definition, the function
\begin{equation}\label{e18}
B(\ix)=e^{-\ix} \sum_{k=1}^\infty \frac{\ix^k}{k!}\, b_k \,. 
\end{equation}
Provided the $b_k$'s grow not too rapidly 
this is an entire function
of $\theta$. In combinatorics, it is usually called the 
exponential generating function of the sequence $\{b_k\}$.
Various methods of extracting asymptotics of sequences from
their generating functions are classically known and widely
used, see for example \cite{V} where such methods are used
to obtain the limit shape of a typical partition under various
measures on the set of partitions. 

A probabilistic way to look at the generating function \eqref{e18} is the
following. If $\theta\ge 0$ then $B(\theta)$ is the expectation
of $b_\eta$ where $\eta\in\{0,1,2,\dots\}$ is a Poisson random
variable with parameter $\theta$. Because $\eta$ has mean $\theta$
and standard deviation $\sqrt\theta$, one expects that 
\begin{equation}\label{e19}
B(n) \approx b_n\,, \quad n\to\infty \,,
\end{equation}
provided the variations of $b_k$ for $|k-n|\le \const\sqrt n$ are small. 
One possible regularity condition on $b_n$ which implies \eqref{e19} 
is monotonicity. In a very general and very 
convenient form, a depoissonization lemma for nonincreasing
nonnegative $b_n$  was established by K.~Johansson in \cite{J}.
We use this lemma in Section \ref{s4} to prove Theorem \ref{t3}.

Another approach to depoissonization is to use a contour integral
\begin{equation}\label{e110}
b_n = \frac{n!}{2\pi i} \int_C \frac{B(z)\, e^{z}}{z^n} \,  \frac{dz}{z} \,,
\end{equation}
where $C$ is any contour around $z=0$.  Suppose, for a moment, that
$b_n$ is constant $b=b_n=B(z)$.  
The function $e^z/z^n=e^{z-n\ln z}$ has a unique critical point $z=n$. If
we choose $|z|=n$ as the contour $C$,  then  only neighborhoods of
size $|z-n|\le \const\sqrt{n}$ contribute to the asymptotics of \eqref{e110}.
Therefore, for general $\{b_n\}$, we still expect that provided  
the overall growth of $B(z)$ is under control  and the
variations of $B(z)$ for $|z-n| \le \const\sqrt n$  are small,
the asymptotically significant contribution to \eqref{e110} will come from 
$z=n$. That is,  we still expect \eqref{e19} to be valid. See, for example,
\cite{JS} for a comprehensive discussion and survey of this approach. 

We use this 
approach to prove Theorem \ref{t2} in Section \ref{s3}. The growth
conditions on $B(z)$ which are suitable in our situation are
spelled out in Lemma \ref{l31}. 

In our case, the functions $B(\ix)$ are combinations of the Bessel
functions. Their asymptotic behavior as $\theta\approx n \to\infty$
can be obtained directly from the classical results on asymptotics 
of Bessel functions which are  discussed, for example, in the fundamental
Watson's treatise \cite{W}. 
These asymptotic formulas for Bessel functions are derived using
the integral representations of Bessel functions and the steepest descent
method. The different behavior of the asymptotics in the bulk $(-2,2)$ of the spectrum,
near the edges $\pm 2$ of the spectrum, and outside of $[-2,2]$  is produced
by the different location of the saddle point in these three cases. 

\subsection{Organization of the paper}
Section \ref{s2} contains the proof of Theorems \ref{t1} and \ref{t1b}
 and also various
formulas for the kernels $\KK$ and $\J$. We also discuss a difference
operator which commutes with $\J$ and its possible applications. 

Section \ref{s3} deals with the
behavior of the Plancherel measure in the bulk of the spectrum;
there we prove Theorem \ref{t2}. Theorem \ref{t3} and a similar
result (Theorem \ref{t4}) for the poissonized measure $M^\ix$
are established in Section \ref{s4}. 

At the end of the paper there is an Appendix, 
where we collected some necessary results about Fredholm determinants,
point processes, and convergence of trace class operators.

\subsection{Acknowledgements}
In many different ways, our work was inspired by the work 
of J.~Baik, P.~Deift, and K.~Johansson, on the one hand, and by the 
work of A.~Vershik and S.~Kerov, on the other. It is our great 
pleasure to thank them for this inspiration and for many 
fruitful discussions.

\section{Correlation functions of the measures $M^\ix$}\label{s2}

\subsection{Proof of Theorem \ref{t1}}\label{s21}

As noted above, Theorem \ref{t1} is a limit case of Theorem 3.3 of \cite{BO}.
That theorem concerns a family $\{M^{(n)}_{zz'}\}$ of probability
measures on partitions of $n$, where $z,z'$ are certain parameters. When the
parameters go to infinity, $M^{(n)}_{zz'}$ tends to the Plancherel measure
$M_n$. Theorem 3.3 in \cite{BO} gives a determinantal 
formula for the correlation functions
of the measure
\begin{equation}\label{e21}
M^\xi_{zz'} = (1-\xi)^t\,\sum_{n=1}^\infty \frac{(t)_n}{n!}\,\xi^n \, M^{(n)}_{zz'} 
\end{equation}
in terms of a certain {\it hypergeometric kernel}.
Here $t=zz'>0$ and $\xi\in(0,1)$ is an additional parameter. 
 As $z,z'\to\infty$
and $\xi=\frac{\ix}{t}\to0$,  the negative binomial distribution in \eqref{e21}
tends to  the Poisson
distribution with parameter $\ix$. In the same limit, the 
hypergeometric kernel becomes the kernel $\KK$ of Theorem \ref{t1}.
The Bessel functions appear as a suitable degeneration of hypergeometric functions.

Recently, these results of \cite{BO} were considerably generalized in \cite{O2},
where it was shown how this type of correlation functions computations
can be done using simple commutations relations in the infinite wedge space. 

For the reader's convenience, we present here a direct and elementary
proof which  uses the same ideas as in \cite{BO} plus an additional technical trick, namely,
differentiation with respect to $\ix$ which kills denominators. This trick
yields an denominator--free integral formula for the kernel $\KK$, see
Proposition \ref{p28}. Our proof here is a verification, not deduction. 
For more conceptual approaches the reader is referred to \cite{BO2,O2}. 

Let $x,y\in\Z+\sh$. Introduce the following kernel $\LL$ 
\begin{equation*}
\LL(x,y;\ix)=\cases 
0\,, & xy>0 \,, \\
\dfrac1{x-y}\dfrac{\ix^{(|x|+|y|)/2}}{\Gamma(|x|+\sh)\, \Gamma(|y|+\sh)} \,, 
& xy <0 \,.
\endcases 
\end{equation*}
We shall 
consider the kernels $\KK$ and $\LL$ as operators in the 
$\ell^2$ space on $\Z+\sh$.  

We recall that simple multiplicative formulas
(for example, the hook formula) are known for the 
number $\dim\la$ in \eqref{e0001}. For our purposes, it is convenient
to rewrite the hook formula in the following determinantal
form. Let $\la=(p_1,\dots,p_d\,|\,q_1,\dots,q_d)$ be the 
Frobenius  coordinates of  $\la$, see Section \ref{s12}. We have 
\begin{equation}\label{e11}
\frac{\dim\la}{|\la|!}=\det\left[\frac1{(p_i+q_j+1)\, p_i!\, q_i!}
\right]_{1\le i,j\le d}\,.
\end{equation}
The following proposition is a straightforward computation using \eqref{e11}. 

\begin{proposition}\label{p21} Let $\la$ be a partition. Then 
\begin{equation}\label{e22}
M^\ix(\la)= e^{-\ix}\, \det \Big[\LL(x_i,x_j;\ix)\Big]_{1\le i,j \le s}\,,
\end{equation}
where $\Fr(\la)=\{x_1,\dots,x_s\}\subset\Z+\sh$ are the 
modified Frobenius coordinates of $\la$. 
\end{proposition}

Let $\Fr_*\left(M^\ix\right)$ be the push-forward of $M^\ix$
under the map $\Fr$. 
Note that the image of  $\Fr$ consists of sets
$X\subset\Z+\sh$ having equally many positive and negative elements.
For other $X\subset\Z+\sh$, the right-hand side of \eqref{e22} can be 
easily seen to vanish. Therefore  $\Fr_*\left(M^\ix\right)$ 
is a determinantal point process  (see the Appendix) 
corresponding to $\LL$, that is, its configuration
probabilities are determinants of the form \eqref{e22}.

\begin{corollary}\label{c22}  $\det(1+\LL)=e^{\ix}$.
\end{corollary}

This follows from the fact that $M^\ix$ is a probability measure. 
This is explained in Propositions
\ref{p51} and \ref{p54} in the Appendix. Note that, in general, 
one needs to check that $\LL$ is a trace class operator. However, because of 
the special form of $\LL$,  it suffices to check a weaker claim -- that $\LL$ is a
Hilbert--Schmidt operator, which is immediate. 

Theorem \ref{t1} now follows from general properties of determinantal point
processes (see Proposition \ref{p55}  in the Appendix) and the following

\begin{proposition}\label{p23}  $\KK=\LL\,(1+\LL)^{-1}$.
\end{proposition}

We shall need three following identities for Bessel functions
which are degeneration of the identities (3.13--15) in \cite{BO}
for the hypergeometric function. The first identity is  due to Lommel 
 (see \cite{W}, Section 3.2 or \cite{HTF}, 7.2.(60))
\begin{equation}\label{e23}
J_\nu(2z)\,J_{1-\nu}(2z) + J_{-\nu}(2z)\, J_{\nu-1}(2z) = \frac{\sin\pi\nu}{\pi\, z} \,.
\end{equation}
The other two identities are the following.  

\begin{lemma}\label{l24} For any $\nu\ne 0,-1,-2,\dots$ and any $z\ne 0$ we have 
\begin{align}
&\sum_{m=0}^\infty  \frac1{m+\nu} \frac{z^m}{m!} \, J_m(2z)  
=\frac{\Gamma(\nu)\, J_\nu(2z)}{z^\nu}  \,,\label{e24}\\
&\sum_{m=0}^\infty  \frac1{m+\nu} \frac{z^m}{m!} \, J_{m+1}(2z)  =
\frac1z-\frac{\Gamma(\nu)\, J_{\nu-1}(2z)}{z^\nu} \label{e25}\,. 
\end{align}
\end{lemma}

\begin{proof} 
Another identity  due to Lommel (see \cite{W}, Section 5.23, or \cite{HTF}, 7.15.(10)) reads
\begin{equation*}
\sum_{m=0}^\infty  \frac{\Gamma(\nu-s+m)}{\Gamma(\nu+m+1)}
\frac{z^m}{m!} \, J_{m+s}(2z)  =
\frac{\Gamma(\nu-s)}{\Gamma(s+1)} \, \frac{J_\nu(2z)}{z^{\nu-s}}   \,.
\end{equation*}
Substituting $s=0$ we get \eqref{e24}. Substituting $s=1$ yields
\begin{equation}\label{e25a}
\sum_{m=0}^\infty  \frac1{(m+\nu)(m+\nu-1)} \frac{z^m}{m!} \, J_{m+1}(2z)  =
\frac{\Gamma(\nu-1)\, J_\nu(2z)}{z^{\nu-1}}  \,.  
\end{equation}
Let $r(\nu,z)$  the difference of the left-hand side and the right-hand side
in \eqref{e25}. Using \eqref{e25a} and the recurrence relation 
\begin{equation}\label{e26}
J_{\nu+1} (2z)-\frac{\nu}{z} J_\nu(2z) + J_{\nu-1} (2z) = 0
\end{equation}
we find that $r(\nu+1,z)=r(\nu,z)$. Hence for any $z$ it is a periodic function
of $\nu$ and it suffices to show that $\lim_{\nu\to\infty} r(\nu,z)=0$. 
Clearly, the left-hand side in \eqref{e25} goes to 0 as $\nu\to\infty$. From
the defining series for $J_\nu$ it is clear that
\begin{equation}\label{e27}
J_\nu(2z)\sim \frac{z^\nu}{\Gamma(\nu+1)} \,, \quad \nu\to\infty \,,  
\end{equation}
which implies that the right-hand side of \eqref{e25} also goes to $0$ as $\nu\to\infty$. 
This concludes the proof. \end{proof}

\begin{proof}[Proof of Proposition] It is convenient to set $z=\sqrt{\ix}$. Since
the operator $1+\LL$ is invertible we have
to check that
\begin{equation*}
\KK+\KK\, \LL - \LL = 0 \,.
\end{equation*}
This is clearly true for $z=0$; therefore, it suffices to check that
\begin{equation}\label{e28}
\dot \KK + \dot \KK \, \LL + \KK\dot \LL -\dot \LL = 0\,,   
\end{equation}
where $\dot \KK = \frac{\partial \KK}{\partial z}$ and $\dot \LL = \frac{\partial \LL}{\partial z}$.  
Using the formulas  
\begin{alignat}{2}\label{e29}
\frac{d}{dz}\, J_x(2z)&=-&&2J_{x+1}(2z)+\frac{x}{z}\, J_{x}(2z) \\
&=&&2J_{x-1}(2z)-\frac{x}{z}\, J_{x}(2z) \nonumber
\end{alignat}
one computes
\begin{equation*}
\dot \KK(x,y)=
\cases
J_{|x|-\h}\, J_{|y|+\h}  + J_{|x|+\h} \, J_{|y|-\h} \,, & xy >0 \,,\\
\sgn(x) \left(J_{|x|-\h}\, J_{|y|-\h}  - J_{|x|+\h} \, J_{|y|+\h} \right) \,, & xy <0\,,
\endcases 
\end{equation*}
where $J_x=J_x(2z)$.  Similarly,
\begin{equation*}
\dot \LL (x,y)=\cases 
0\,, & xy>0 \,, \\
\sgn(x)\, \dfrac{z^{|x|+|y|-1}}{\Gamma(|x|+\sh)\, \Gamma(|y|+\sh)} \,, 
& xy <0 \,. 
\endcases 
\end{equation*}
Now the verification of \eqref{e28} becomes a straightforward application of the
formulas \eqref{e24} and \eqref{e25}, except for the occurrence of the singularity $\nu\in\Z_{\le 0}$
in those formulas. This singularity is resolved using \eqref{e23}. 
This concludes the proof of Proposition \ref{p23} and Theorem \ref{t1}. \end{proof}

\subsection{Proof of Theorem \ref{t1b}}\label{s21b}

Recall that by construction 
$$
\Fr(\la)=\left(\sD(\la)+\sh\right)\triangle \left(\Z_{\le 0} -\sh\right)\,.
$$
Let us check that this and Proposition \ref{pAc} implies Theorem \ref{t1b}.
In Proposition \ref{pAc} we substitute 
$$
\X=\Z+\sh\,, \quad Z=\Z_{\le 0}-\sh\,, \quad K=\KK\,.
$$
By definition, set
$$
\ep(x)=\sgn(x)^{x+1/2}\,, \quad x\in\Z+\sh\,.
$$
We have the following

\begin{lemma}\label{l0001} $\KK^\tri(x,y) = \ep(x)\, \ep(y) \, \J(x-\sh,y-\sh)$
\end{lemma}

It is clear that since the $\ep$-factors cancel out of all determinantal
formulas, this lemma and Proposition \ref{pAc} establish  the equivalence
of Theorems \ref{t1} and \ref{t1b}.

\begin{proof}[Proof of lemma]
Using the relation 
$$
J_{-n}=(-1)^n J_n
$$
and the definition of $\KK$ one computes
\begin{equation}\label{e210}
\KK(x,y) = \sgn(x)\, \ep(x)\, \ep(y) \, \J(x-\sh,y-\sh)\,, \quad x\ne y\,.
\end{equation}
Clearly, the relation \eqref{e210} remains valid for $x=y>0$. It remains
to consider the case $x=y<0$. In this case we have to show that
$$
1-\KK(x,x)=\J(x-\sh,y-\sh)\,, \quad x\in\Z_{\le 0}-\sh \,.
$$
Rewrite it as
\begin{equation}\label{e004}
1-\J(k,k)=\J(-k-1,-k-1)\,, \quad k=-x-\sh\in\Z_{\ge 0}\,.
\end{equation}
By \eqref{e213} this is equivalent to
\begin{multline*}
1-\sum_{m=0}^\infty (-1)^m\,\frac{(2k+m+2)_m}{\Ga(k+m+2)\Ga(k+m+2)}\,
\frac{\theta^{k+m+1}}{m!}\\
=\sum_{n=0}^\infty (-1)^n\,\frac{(-2k+n)_n}{\Ga(-k+n+1)\Ga(-k+n+1)}\,
\frac{\theta^{-k+n}}{n!}\,.
\end{multline*}
Examine the right--hand side. The terms with $n=0,\dots,k-1$ vanish
because then $1/\Ga(-k+n+1)=0$. The term with $n=k$ is equal to 1,
which corresponds to 1 in the left--hand side. Next, the terms with
$n=k+1,\dots,2k$ vanish because for these values of $n$, the expression
$(-2k+n)_n$ vanishes. Finally, for $n\ge 2k+1$, set $n=2k+1+m$. Then the
$n$th term in the second sum is equal to minus the $m$th term in the
first sum. Indeed, this follows from the trivial relation 
$$
-(-1)^m\,\frac{(2k+m+2)_m}{m!}=(-1)^n\,\frac{(-2k+n)_n}{n!}\,,
\qquad n=2k+1+m.
$$
This concludes the proof.
\end{proof}

\subsection{Various formulas for the kernel $\J$}\label{s22}

Recall that since $J_x$ is an entire 
function of $x$, the function $\J(x,y)$ is entire in $x$ and $y$. 
We shall now obtain several denominator--free formulas
for the kernel $\J$. 

\begin{proposition}\label{p27} 
\begin{equation}\label{e213}
\J(x,y;\ix)=\sum_{m=0}^\infty(-1)^m\,
\frac{(x+y+m+2)_m}{\Ga(x+m+2)\Ga(y+m+2)}\,
\frac{\ix^{\frac{x+y}2+m+1}}{m!}\,.
\end{equation}
\end{proposition}

\begin{proof} Straightforward computation using a formula due
to Nielsen, see Section 5.41 of \cite{W} or 
\cite{HTF}, formula 7.2.(48) . \end{proof}

\begin{proposition}\label{p28} Suppose $x+y>-2$. Then 
\begin{equation*}
\J(x,y;\ix)=\frac12\, \int_0^{2\sqrt\ix} (J_x(z)\, J_{y+1}(z) + J_{x+1}(z)\, J_{y}(z))\, dz
\end{equation*}
\end{proposition}
\begin{proof} Follows from a computation done in the proof of Proposition \ref{p23}
\begin{equation*}
\frac{\partial}{\partial\ix} \, \J(x,y;\ix) = \frac{1}{2\sqrt{\ix}} \,
(J_x \, J_{y+1} + J_{x+1}\, J_{y}) \,,
\quad J_x=J_x(2\sqrt{\ix})\,, 
\end{equation*}
and the following corollary of \eqref{e213}
$$
\J(x,y;0)=0\,, \quad x+y>-2\,.
$$
\end{proof}

\begin{remark} 
Observe that  by Proposition \ref{p28} the operator $\frac{\partial\J}{\partial\ix}$
is a sum of two operators of rank 1. 
\end{remark}

\begin{proposition}\label{p29} 
\begin{equation}\label{e214}
\J(x,y;\ix)=\sum_{s=1}^\infty
J_{x+s}\, J_{y+s}\,, \qquad J_x=J_x(2\sqrt\ix).  
\end{equation}
\end{proposition}

\begin{proof} Our argument is similar to an argument due to
Tracy and Widom, see the proof of the formula (4.6)
in \cite{TW}.  The recurrence relation \eqref{e26} implies that 
\begin{equation}\label{e215}
\J(x+1,y+1)-\J(x,y)=-J_{x+1}\, J_{y+1}
\end{equation}
Consequently, the difference between the left-hand side and 
the right-hand side  of \eqref{e214} is a
function which depends only on $x-y$. Let $x$ and $y$ go to infinity 
in such a way that $x-y$ remains fixed. Because of 
the asymptotics \eqref{e27} both sides in \eqref{e214}
tend to zero and, hence, the difference actually is 0. \end{proof}

In the same way as in \cite{TW} this results in  the following 

\begin{corollary}\label{c210} For any
$a\in\Z$, the restriction of the 
kernel $\J$ to the subset $\{a,a+1,a+2,\dots\}\subset\Z$
determines a nonnegative trace class operator in the $\ell^2$
space on that subset.
\end{corollary}

\begin{proof} By Proposition \ref{p29}, the restriction of $\J$ on
$\{a,a+1,a+2,\dots\}$ is the square of the kernel 
$(x,y)\mapsto J_{x+y+1-a}(2\sqrt\ix)$. Since the latter kernel is
real and symmetric, the kernel  $\J$ is nonnegative.
Hence, it remains to prove that its trace is finite. Again, by
Proposition \ref{p29}, this trace is equal to 
\begin{equation*}
\sum_{s=1}^\infty s \, (J_{a+s+1}(2\sqrt\ix))^2.
\end{equation*}
This sum is clearly finite by \eqref{e27}. 
\end{proof}

\begin{remark}
The kernel $\J$ resembles
a Christoffel--Darboux kernel and, in fact, 
the operator in $\ell^2(\Z)$ defined by the
kernel $\J$ is an Hermitian projection operator.
Recall that $\KK=\LL(1+\LL)^{-1}$, where $\LL$ 
is of the form
$$
\LL=\begin{bmatrix} 0 & A \\ -A^* & 0 
\end{bmatrix}
$$
On can prove that this together with Lemma \ref{l0001} implies
that  $\J$ is an Hermitian projection kernel.
However, in contrast to a Christoffel--Darboux kernel, it
 projects to an infinite--dimensional subspace. 
\end{remark}

\subsection{Commuting difference operator}

Consider the difference operators $\De$ and $\nabla$ on the lattice $\Z$,
$$
(\De f)(k)=f(k+1)-f(k)\,, \qquad (\nabla f)(k)=f(k)-f(k-1)\,. 
$$
Note that $\nabla=-\De^*$ as operators on $\ell^2(\Z)$. Consider the following
 second order difference
Sturm--Liouville operator 
\begin{equation}\label{e0002}
D=\De\circ\al\circ\nabla +\be\,,
\end{equation}
where $\al$ and $\be$ are operators of multiplication by certain
functions $\al(k)$, $\be(k)$. The operator \eqref{e0002} is
self--adjoint in $\ell^2(\Z)$. A straightforward computation shows that
\begin{multline}\label{e0003}
\big[D f\big](k)=(-\al(k+1)-\al(k)+\be(k))f(k)+\\
\al(k)f(k-1)+\al(k+1)f(k+1)\,.  
\end{multline}
It follows that if $\al(s)=0$ for a certain $s\in\Z$ then the space
of functions $f(k)$ vanishing for $k<s$ is invariant under $D$.

\begin{proposition} Let $[\J]_s$ denote the operator in
$\ell^2(\{s,s+1,\dots\})$ obtained by restricting the kernel
$\J$ to $\{s,s+1,\dots\}$. Then the difference
Sturm--Liouville operator \eqref{e0002} commutes with $[\J]_s$ provided
$$
\al(k)=k-s, \qquad
\be(k)=-\,\frac{k(k+1-s-2\sqrt\ix)}{\sqrt\ix}+\textup{const} \,.
$$
\end{proposition}

\begin{proof} Since $[\J]_s$ is the square of the operator with the
kernel $J_{k+l+1-s}$, it suffices to check that the latter
operator commutes with $D$, with the above choice of  $\al$ and $\be$.
But this is readily checked using \eqref{e0003}. 
\end{proof}

This proposition is a counterpart of a known fact about the Airy kernel, see
\cite{TW}. Moreover, in the scaling limit when  $\ix\to\infty$ and 
$$
k=2\sqrt\ix+x\,\ix^{1/6}, \qquad s=2\sqrt\ix+\varsigma\,\ix^{1/6}, 
$$
the difference operator $D$
becomes, for a suitable choice of the constant, the differential operator 
$$
\frac d{dx}\circ(x-\varsigma)\circ \frac d{dx} -x(x-\varsigma), 
$$
which commutes to the Airy operator restricted to $(\varsigma,+\infty)$.  The
above differential 
operator  is exactly that of Tracy and Widom \cite{TW}.

\begin{remark}
Presumably, this commuting difference operator can be used to
obtain, as was done in \cite{TW} for the Airy kernel, 
asymptotic formulas for the eigenvalues of $[\J]_s$, where
$s=2\sqrt\ix+\varsigma\,\ix^{1/6}$ and $\varsigma\ll 0$. 
Such asymptotic formulas may be very useful if one wishes
to refine Theorem \ref{t3} and to establish
convergence of moments in addition to convergence of 
distribution functions. For individual distributions of
$\la_1$ and $\la_2$ the convergence of moments was 
obtained, by other methods, in \cite{BDJ1,BDJ2}. 
\end{remark}

\section{Correlation functions in the bulk of the spectrum}\label{s3}

\subsection{Proof of Theorem \ref{t2}}

We refer the reader to Section \ref{s13} of the Introduction for
the definition of a regular sequence $X(n)\subset\Z$ and the
statement of Theorem \ref{t2}. Also, in this section, we shall be working in the bulk of the spectrum,
that is, we shall assume that all numbers $a_i$ defined in \eqref{e14}
lie inside $(-2,2)$. The edges $\pm2$ of the spectrum and its exterior will be treated  in
the next section. 
 
In our proof, we shall follow the strategy explained in Section \ref{s15}.
Namely, in order to compute the limit of $\brh(n,X(n))$ we shall
use the contour integral
\begin{equation*}
\brh(n,X(n)) = \frac{n!}{2\pi i} \int_{|\ix|=n}  \brh^\ix(X(n))\, 
\frac{e^\ix}{\ix^{n+1}} \,  d\ix \,, 
\end{equation*}
compute the asymptotics of $\brh^\ix$ for $\ix\approx n$, and
estimate $|\brh^\ix|$ away from $\ix=n$.  Both tasks will be 
accomplished using classical results about the Bessel functions. 

We start our proof  with the following
lemma which  formalizes the above  informal depoissonization
argument.  The hypothesis of this lemma
is very far from optimal, but it is sufficient for our purposes.
For the rest of this section, we fix a number $0<\al<1/4$ which
shall play an auxiliary role.   

\begin{lemma}\label{l31} Let $\{f_n\}$ be a sequence of entire functions
\begin{equation*}
f_n(z)=e^{-z} \sum_{k\ge 0} \frac{f_{nk}}{k!} \, z^k \,,
\quad n =1,2,\dots \,, 
\end{equation*}
and suppose that there exist  such constants  $f_\infty$ and $\g$ that
\begin{align} 
&\max_{|z|= n} \left|f_n(z)\right|  
= O\left(e^{\g\,\sqrt{n}}\right)  \label{A}
\\
&
\max_{|z/n-1|\le n^{-\al}} 
\left|f_n(z)-f_\infty\right| e^{-\g|z-n|/\sqrt{n}} = o(1) \label{B}\,,
\end{align} 
as $n\to\infty$.  Then 
\begin{equation*}
\lim_{n\to\infty} f_{nn} =  f_\infty\,.
\end{equation*}
\end{lemma}

\begin{proof} By replacing $f_n(z)$ by $f_n(z)-f_\infty$, we
may assume that $f_\infty=0$. 
By Cauchy and Stirling formulas, we have
\begin{equation*}
f_{nn} =(1+o(1))\,  \sqrt{\frac{n}{2\pi}} \, \int_{|\zeta|=1} 
\frac{f_n(n\zeta)\, e^{n(\zeta-1)}}{\zeta^n} \, \frac{d\zeta}{i\zeta}  \,.
\end{equation*}
Choose some large $C>0$  and split the 
circle $|\zeta|=1$ into 2  parts as follows:
\begin{equation*}
S_1=\left\{ \frac{C}{n^{1/4}}\le |\zeta-1|\right\}\,, \quad
S_2=\left\{ \frac{C}{n^{1/4}}\ge|\zeta-1|\right\} \,.
\end{equation*}
The inequality \eqref{A} and the equality  
\begin{equation*}
\left|e^{n(\zeta-1)}\right|=e^{-n|\zeta-1|^2/2} \,.
\end{equation*}
imply that  the integral $\int_{S_1}$ decays exponentially provided $C$ is large enough. 
On $S_2$,  the inequality \eqref{B} applies for sufficiently large $n$ and gives 
\begin{equation*}
\max_{z\in S_2} \left|f_n(n\zeta)\right| e^{-\g\sqrt{n}|\zeta-1|} = o(1) \,.
\end{equation*}
Therefore, the the integral $\int_{S_2}$ is  $o(\,\,)$ of the following integral 
\begin{equation*}
\sqrt{n} \int_{|\zeta|=1} \frac{d\zeta}{i\zeta} \,
\exp\left(-n\frac{|\zeta-1|^2}2 + \g\sqrt{n}|\zeta-1|\right) 
\sim \int_{-\infty}^\infty e^{-s^2/2 +\gamma |s|} \, ds \,.
\end{equation*}
Hence,  $\int_{S_2}=o(1)$ and the lemma follows.  
\end{proof}

\begin{definition} Denote by $\aF$ the algebra
(with respect to termwise addition and multiplication)
 of sequences $\{f_n(z)\}$ which satisfy the 
properties \eqref{A} and \eqref{B} for some, depending on the
sequence,  constants $f_\infty$ and $\g$. Introduce
the following map
\begin{equation*}
\Lim:\aF\to \C\,, \quad \{f_n(z)\} \mapsto f_\infty \,,
\end{equation*}
which is clearly  a homomorphism.  
\end{definition} 

\begin{remark}
Note that we do not require $f_n(z)$ to be entire. Indeed, the
kernel $\J$ may have a square root branching, see
the formula \eqref{e213}.
\end{remark} 

By Theorem \ref{t1b},  the correlation 
functions $\brh^\ix$ belong to the algebra generated by
sequences of the form
\begin{equation*}
\{f_n(z)\}=\left\{\J(x_n,y_n;z)\right\}\,,   
\end{equation*}
where the sequence  $X=X(n)=\{x_n,y_n\}\subset \Z$ is regular which,
we recall, means that the limits 
\begin{equation*}
a=\lim_{n\to\infty} \frac{x_n}{\sqrt{n}}\,, \quad d=\lim_{n\to\infty}  (x_n-y_n) 
\end{equation*}
exist, finite or infinite. 
Therefore, we first consider such sequences. 

\begin{proposition}\label{p34} 
If $X=\{x_n,y_n\}\subset\Z$ is regular then 
\begin{equation*}
\left\{\J(x_n,y_n;z)\right\}\in\aF\,, \quad 
\Lim\left(\left\{\J(x_n,y_n;z)\right\}\right)= \bS(d,a) \,.
\end{equation*}
\end{proposition} 

In the proof of this proposition it will be convenient to allow $X\subset \C$.
For complex sequences $X$ we shall require $a\in\R$; the number $d\in\C$
may be arbitrary. 

\begin{lemma}\label{l35}  
Suppose that a sequence $X\subset\C$
is as above and, additionally, suppose that $\Im x_n$, $\Im y_n$ are 
bounded and $d\ne 0$. Then the sequence $\left\{\J(x_n,y_n;z)\right\}$ 
satisfies \eqref{B}  with $f_\infty = \bS(d,a)$ and certain $\g$. 
\end{lemma}

\begin{proof}[Proof of Lemma]
We shall use Debye's asymptotic formulas for  Bessel functions of complex
order and large complex argument, see, for example,  Section 8.6 in \cite{W}.
Introduce the following function
\begin{equation*}
F(x,z)=z^{1/4} \, J_x(2\sqrt{z}) \,.
\end{equation*}
The formula \eqref{e212} can be rewritten as follows 
\begin{equation}\label{e31}
\J(x,y;z)=\frac{F(x,z)\, F(y+1,z) - F(x+1,z) \,F(y,z)}{x-y} \,. 
\end{equation}
The asymptotic formulas for Bessel functions imply that 
\begin{equation}\label{e32}
F(x,z)=\frac{\cos
\left(\sqrt{z} \, G(u) +\frac{\pi}4\right)}{H(u)^{1/2}}
\left(1+O\left(z^{-1/2}\right)\right) \,, \quad u=\frac{x}{\sqrt{z}} \,,
\end{equation}
where
\begin{equation*}
G(u)=\frac{\pi}{2}\left(u-\Omega(u)\right)\,,
\quad H(u)=\frac{\pi}2 \, \sqrt{4-u^2} \,,
\end{equation*}
provided that  $z\to\infty$ in such a way that $u$ stays in some neighborhood
of $(-2,2)$; the precise form of this neighborhood can be seen in Figure 22
in Section 8.61 of \cite{W}.  Because we assume that 
\begin{equation*}
\lim_{n\to\infty} \frac{x_n}{\sqrt{n}}\,,\,\lim_{n\to\infty} \frac{y_n}{\sqrt{n}} \in (-2,2) \,,
\end{equation*}
and because $|z/n-1|<n^{-\al}$, 
 the ratios ${x_n}/{\sqrt{z}}$, ${y_n}/{\sqrt{z}}$ stay 
close to $(-2,2)$. For future reference, we also point
out that the constant in $O\left(z^{-1/2}\right)$ in \eqref{e32} is uniform in 
$u$ provided  $u$ is bounded away from the endpoints $\pm2$. 

First  we estimate $\Im \left(\sqrt{z} \, G(u)\right) $. The function
$G$ clearly takes real values on the real line. From the obvious 
estimate 
\begin{equation*}
\left| \Im \left(\sqrt{z} \, G(u) \right) \right|
\le
\left| \Im \left(\sqrt{n} \, G(x/\sqrt{n}) \right) \right|
+ \left| \sqrt{z}  \, G(x/\sqrt{z}) -
\sqrt{n}  \, G(x/\sqrt{n}) \right| 
\end{equation*}
and the boundedness of  $G$, $G'$, and  $|\Im x|$ we obtain an estimate of the form  
\begin{equation}\label{e33}
\max_{|z/n-1|\le n^{-\al}} |F(x;z)| e^{-\const\, |z-n|/\sqrt{n}}  
= O(1)\,. 
\end{equation}

 If $d=\infty$ then 
because of the denominator in \eqref{e31} the estimate \eqref{e33} implies that 
\begin{equation*}
\J(x_n,y_n;z)=o\left(e^{\const\, |z-n|/\sqrt{n}}\right)\,.
\end{equation*}
Since $\bS(\infty,a)=0$, it follows that  in this case the 
lemma is established. 

Assume, therefore, that $d$ is finite. 
Observe that for any bounded increment $\De x$ we  have
\begin{multline}\label{e34}
F(x+\De x,z)=\frac{\cos
\left(\sqrt{z} \, G(u) +G'(u)\,\De x+ \frac{\pi}4\right)}
{H(u)^{1/2}}\\
 +O\left( \frac{(\De x)^2}{\sqrt z}\, e^{\const\, |z-n|/\sqrt{n}}\right)\,,
\end{multline}
and, in particular, the last term is $o\left(e^{\const\, |z-n|/\sqrt{n}}\right)$. 
Using the trigonometric identity 
\begin{equation*}
 \cos\left(A\right)\cos\left(B+C\right)-\cos\left(A+C\right)\cos\left(B\right)=
\sin\left(C\right)\sin\left(A-B\right)\,,  
\end{equation*}
and observing that
\begin{equation*}
G'(u)=\arccos(u/2) \,, 
\quad \sin(G'(u))=\frac{\sqrt{4-u^2}}2=\frac{H(u)}{\pi}\,,
\end{equation*}
we compute 
\begin{multline*}
F(x_n;z)\, F(y_n+1;z) -F(x_n+1;z) \,F(y_n;z) = \\
\frac1{\pi} \sin\left(\arccos\left(\frac{x_n}{2\sqrt{z}}\right)\, (x_n-y_n)\right) + 
 o\left(e^{\const\, |z-n|/\sqrt{n}}\right)\,.
\end{multline*} 
Since, by  hypothesis,
\begin{equation*}
\frac{x_n}{\sqrt{z}} \to a\,, \quad  (x_n-y_n) \to d \,, 
\end{equation*}
and $d\ne 0$, the lemma follows. \end{proof}
 
\begin{remark}
Below  we shall need this lemma for a variable  sequence $X=\{x_n,y_n\}$.
Therefore, let us spell out explicitly under what conditions on $X$  the
estimates in Lemma \ref{l35} remain uniform. We need the sequences
$\frac{x_n}{\sqrt{n}}$ and $\frac{y_n}{\sqrt{n}}$ to converge uniformly; then, in particular,
the ratios $\frac{x_n}{\sqrt{n}}$ and $\frac{y_n}{\sqrt{n}}$ are uniformly bounded
away from $\pm 2$.  Also, we need $\Im x_n$ and $\Im y_n$ to be
uniformly bounded. Finally, we need $|d|$ to be uniformly bounded from below.
\end{remark}

\begin{proof}[Proof of Proposition] First, we check the condition \eqref{B}. In the case
$d\ne 0$ this was done in the previous lemma.  Suppose, therefore,
that $\{x_n\}$ is a regular sequence in $\Z_{\ge 0}$ and
consider the asymptotics of $\J(x_n,x_n;z)$.

Because the function $\J(x,y;z)$ is an entire function of $x$ and $y$ we have
\begin{equation}\label{e35}
\J(x,x;z) = \frac1{2\pi} \int_0^{2\pi} \J\left(x,x+r e^{i t};z\right) \, dt \,,  
\end{equation}
where $r$ is arbitrary; we shall take $r$ to be some small but fixed
number. From the previous lemma we know that 
\begin{equation*}
\J\left(x,x+r e^{i t};z\right)=
\frac1{\pi r e^{it}}
\sin\left(\omega\left(\frac{x}{\sqrt{z}}\right)\, re^{it}\right)
 + o\left(e^{\const\, |z-n|/\sqrt{n}}\right) \,.
\end{equation*}
{}From the above remark it follows that this estimate is uniform in $t$.
This implies the property \eqref{B} for 
$\J(x_n,x_n;z)$. 

To prove the estimate \eqref{A} we use 
Schl\"afli's integral representation (see Section 6.21 in \cite{W})
\begin{multline}\label{e36}
J_x(2\sqrt{z})=\frac1\pi \int_0^\pi \cos
\left(x t - 2 \sqrt{z} \, \sin t\right) \, dt - \\
\frac{\sin \pi x}{\pi} \int_0^\infty e^{-x t -2\sqrt{z} \, \sinh t} \, dt \,, 
\end{multline}
which is valid for $|\arg z|<\pi$ and even for $\arg z=\pm \pi$ 
provided $\Re x>0$ or $x\in\Z$. 

If $x\in\Z$ then the second summand in \eqref{e36}
vanishes and and the first is $O\left(e^{\const|z|^{1/2}}\right)$ uniformly
in $x\in\Z$. This implies the estimate \eqref{A} 
provided $d\ne 0$. 

It remains, therefore, to check \eqref{A} for $\J(x_n,x_n;z)$ where
$\{x_n\}\in\Z$ is a regular sequence. Again, we use 
\eqref{e35}. Observe, that since $\Re\sqrt{z}\ge 0$ 
the second summand  in \eqref{e36} is uniformly small
provided $\Im x$ is bounded from above and $\Re x$ is bounded
from below. Therefore, \eqref{e35} produces the \eqref{A} estimate 
for $x_n\ge 1$. For $x_n\le 0$ we use the relation \eqref{e004} and 
the reccurence \eqref{e215} to obtain the estimate.  
\end{proof}

\begin{proof}[Proof of Theorem \ref{t2}]
Let $X(n)$ be a regular sequence and let the numbers $a_i$ and
$d_{ij}$ be defined by \eqref{e14}, \eqref{e15}.  We shall assume that $|a_i|<2$ for
all $i$. The validity of the theorem in the case when $|a_i|\ge 2$ for some $i$
will be obvious form the results of the next section. 

We have
\begin{align}
\brh^\ix(X(n))&=e^{-\ix}\, \sum_{k=0}^\infty \brh(k,X(n))\, \frac{\ix^k}{k!}  \label{e37}\\
&=\det \Big[\J(x_i(n),x_j(n))\Big]_{1\le i,j \le s}\,.  \label{e38}
\end{align}
where the first line is the definition of $\brh^\ix$ and the second
is Theorem \ref{t1b}.  From \eqref{e37} it is obvious that $\brh^\ix$ is
entire. Therefore, we can apply Lemma \ref{l31} to it. It is
clear that Lemma \ref{l31}, together with Proposition \ref{p34}, implies
Theorem \ref{t2}. The factorization \eqref{e1t2} follows from 
the vanishing  $\bS(\infty,a)=0$. 
\end{proof}

\subsection{Asymptotics of $\rho(n,X)$}\label{s32}

Recall that the correlation functions $\rho(n,X)$ were
defined by
\begin{equation*}
\rho(n,X)= M_n\left(\left\{\la\,|\, X\subset \Fr(\la)\right\} 
\right)\,, \quad X\subset\Z+\sh\,. 
\end{equation*}
The asymptotics of these correlation functions can be easily
obtained from Theorem \ref{t2} by complementation, see
Sections \ref{sA3} and \ref{s21b}, and the result is the 
following.

Let $X(n)\subset \Z+\sh$ be a regular sequence. If it
splits, then the limit $\lim_{n\to\infty} \rho(n,X(n))$ factors
as in \eqref{e1t2}. Suppose therefore, that $X(n)$
is nonsplit. Here one has to distinguish two cases.
If $X(n)\subset\Z_{\ge 0}+\sh$ or $X(n)\subset\Z_{\le 0}-\sh$
then we shall say that this sequence is \emph{off-diagonal}.
Geometrically, it means that $X(n)$ corresponds to
modified Frobenius coordinates of only one kind:
either the row ones or the column ones. For off-diagonal
sequences we obtain from Theorem \ref{t2} by complementation
that
\begin{equation*}
\lim_{n\to\infty} \rho\left(n,X(n)\right)
=
\det \Big[\,\bS(d_{ij},|a|)
\Big]_{1\le i,j \le s}\,,
\end{equation*}
where $\bS$ is the discrete sine kernel and $a=a_1=a_2=\dots$.

If $X(n)$ is nonsplit and \emph{diagonal}, that is, if it
is nonsplit and includes both positive and negative numbers,
then one has to assume additionally that the number of
positive and negative elements of $X(n)$ stabilizes for 
sufficiently large $n$. In this case the limit correlations
are given by the kernel
\begin{equation}\label{edefD}
\DD(x,y)=\cases
\bS\left(x-y,0\right) \,, & xy>0\,, \\
\displaystyle \frac{\cos\left(\frac\pi 2(x+y)\right)}{\pi(x-y)} \,,  &xy<0 \,.
\endcases
\end{equation}
Remark that this kernel \emph{is not} translation invariant. Note,
however, that
$$
\DD(x+1,y+1)=\sgn(xy)\, \DD(x,y)\,,
$$
provided $x$ and $x+1$ have the same sign and similarly for $y$. 
Therefore, if the subsets $X\subset\Z+\sh$ and $X+m$, $m\in\Z$,
have the same number of positive and negative elements then
$$
\det \Big[ \DD(x_i,x_j)\Big]_{x_i\in X} =
\det \Big[ \DD(x_i+m,x_j+m)\Big]_{x_i\in X} \,.
$$

\section{Edge of the spectrum: convergence to the Airy ensemble}\label{s4}

\subsection{Results and strategy of proof}\label{s41}

In this section we prove Theorem \ref{t3} which was stated in Section \ref{s14}
of the Introduction. We refer the reader to Section \ref{s14} for a discussion
of the relation between Theorem \ref{t3} and the results obtained in \cite{BDJ1,BDJ2,O}.  

Recall that the Airy kernel was defined as follows 
\begin{equation*}
\A(x,y)=\frac{A(x)A'(y)-A'(x)A(y)}{x-y}.
\end{equation*}
where $A(x)$ is the Airy function \eqref{e17}. 
The {\it Airy ensemble} is, by definition,  a random point process on $\R$, whose 
correlation functions are given by 
\begin{equation*}
\rho^{\textup{Airy}}_k(x_1,\dots,x_k)= \det \Big[\,\A(x_i,x_j)
\Big]_{1\le i,j \le k}\,.
\end{equation*}
This ensemble was studied in \cite{TW}. 
We denote by $\zeta_1 > \zeta_2 > \dots $ a random configuration
of the Airy ensemble. 
Theorem \ref{t3} says that after a proper scaling and normalization, the rows
$\la_1,\la_2,\dots$ of a Plancherel random partition $\la$ 
converge in joint distribution 
to the Airy ensemble.  Namely, the following random variables $\widetilde{\la}$
\begin{equation*}
\widetilde{\la}=
\left(\widetilde{\la}_1 \ge \widetilde{\la}_2 \ge \dots
\right)\,,
\quad 
\widetilde{\la}_i =
n^{1/3} \, \left(
\frac{\la_i}{n^{1/2}}-2 \right) \,,
\end{equation*}
converge, in joint distribution, to the Airy ensemble as $n\to\infty$. 

In the proof of Theorem \ref{t3}, we shall follow the strategy explained in Section \ref{s15}
of the Introduction. First, we shall prove that under 
the poissonized measure $M^\ix$ on the set of partitions $\la$, the random variables
$\widetilde{\la}$ 
 converge, in joint distribution,  
to the Airy ensemble as $\ix\approx n\to\infty$.  This result is stated below as
Theorem \ref{t4}. From this, using certain monotonicity and 
Lemma \ref{l48} which is  due to K.~Johansson, we shall conclude 
that the same is true for the measures $M_n$ as $n\to\infty$. 

The proof of Theorem \ref{t4} will be based on the analysis of
 the behavior of the correlation functions
 of $M^\ix$, $\ix\approx n \to \infty$, near the point $2\sqrt{n}$.  From 
the expression for correlation functions of $M^\ix$ given in Theorem \ref{t1}
it is clear that this amounts to the study of  the asymptotics of $J_{2\sqrt{n}}(2\sqrt\ix)$
when $\ix\approx n \to \infty$.  This asymptotics is classically known
and from it we shall derive the following 

\begin{proposition}\label{p41} Set $r=\sqrt\ix$. We have 
\begin{equation*}
r^{\frac 13}\J\left(2r+xr^{\frac 13},2r+yr^{\frac 13},{r^2} \right)\to \A(x,y), 
\quad r\to+\infty\,,
\end{equation*}
uniformly in $x$ and $y$ on compact sets of $\R$.
\end{proposition} 

The prefactor $r^{\frac 13}$ corresponds to the fact that we change the local
scale  near $2r$  to get  non-vanishing limit correlations.  

Using this and verifying certain tail estimates we obtain the following 

\begin{theorem}\label{t4}
For any fixed $m=1,2,\dots$ and any $a_1,\dots,a_m\in\R$ we have 
\begin{multline}\label{e401} 
\lim_{\ix\to+\infty}
M^\ix\left(\left\{\la\,\left|\,
\frac{\la_i-2\sqrt{\ix}}{\ix^{\frac 16}}\,< a_i\,, \  1\le i\le 
m\right.\right\}\right)=\\
\Prob\{\zeta_i< a_i\,,\  1\le i\le m\}\,,
\end{multline}
where $\zeta_1>\zeta_2>\dots$ is the Airy ensemble. 
\end{theorem} 

Observe that the limit behavior of $\widetilde{\la}$ is, obviously, identical
with the limit behavior of similarly scaled $1$st, $2$nd, an so on maximal
Frobenius coordinates. 

Proofs of Proposition \ref{p41} and Theorem \ref{t4} are given Section \ref{s42}.
In Section \ref{s43}, using a depoissonization argument based on Lemma \ref{l48} 
we deduce Theorem \ref{t3}. 

\begin{remark}  We consider the behavior of any number of
initial rows of $\la$, where $\la$
is a Plancherel random partition. By symmetry, same results 
describe the behavior of any number of initial columns of $\la$. 
\end{remark}

\subsection{Proof of Theorem \ref{t4}}\label{s42}

Suppose we have a point process on $\R$ with
determinantal correlation functions 
\begin{equation*}
\rho_k(x_1,\dots,x_k)=\det[K(x_i,x_j)]_{1\le i,j \le k}\,,
\end{equation*}
for some kernel $K(x,y)$.  Let $I$ be a possibly infinite interval $I\subset\R$. 
By  $[K]_I$ we denote the operator in $L^2(I,dx)$
obtained  by restricting the kernel on $I\times I$. Assume $[K]_I$ is
a trace class operator. Then the intersection of the random configuration
$X$ with $I$ is finite almost surely  and
\begin{equation*}
\Prob\{|X\cap I|=N\}=\left.\frac{d^N}{dz^N}\det\Big(1-z[K]_I\Big)\right|_{z=1}\,.
\end{equation*}
In particular, the probability that $X\cap I$ is empty is equal to
\begin{equation*}
\Prob\{X\cap I=\emptyset\}=\det\Big(1-[K]_I\Big)\,.
\end{equation*}

More generally, if $I_1,\dots,I_m$ is a finite family of pairwise
nonintersecting intervals such that the operators
$[K]_{I_1},\dots,[K]_{I_m}$ are trace class then
\begin{multline}\label{e41} 
\Prob\{|X\cap I_1|=N_1,\dots,|X\cap I_m|=N_m \}\\
=\left.\frac{\partial^{N_1+\dots+N_m}}
{\partial z_1^{N_1}\dots\partial z_m^{N_m}}
\det\Big(1-z_1[K]_{I_1}-\dots-z_m[K]_{I_m}\Big)\right|_{z_1=\dots=z_m=1}.
\end{multline}
Here operators $\{[K]_{I_i}\}$ are considered to be acting in the same Hilbert 
space, for example, in $L^2(I_1\sqcup I_2,\sqcup\dots\sqcup I_m,dx)$.

In case of intersecting intervals $I_1,\dots,I_m$, the probabilities
$$
\Prob\{|X\cap I_1|=N_1,\dots,|X\cap I_m|=N_m \}
$$ 
are finite
linear combinations of expressions of the form \eqref{e41}. Therefore,
in order to show the convergence in distribution of point processes with determinantal 
correlation functions, it suffices to show the convergence of expressions of 
the form \eqref{e41}.

The formula \eqref{e41} is discussed, for example, in \cite{TW2}. It 
remains valid for processes on  a lattice such as $\Z$
in  which case the kernel $K$ should be an operator on $\ell^2(\Z)$. 

As verified, for example, in Proposition \ref{p58} in the Appendix, the right-hand
side of \eqref{e41} is continuous in each $[K]_{I_i}$ with respect to the trace norm.
We shall show that after a suitable embedding of $\ell^2(\Z)$ inside
$L^2(\R)$ the kernel $\J(x,y;\ix)$ converges to the Airy kernel $\A(x,y)$
as $\ix\to\infty$.  

Namely, we shall consider a family of embeddings
$\ell^2(\Z)\to L^2(\R)$,  indexed by a positive number $r>0$, which are
defined by
\begin{equation}\label{e402}
\ell^2(\Z)\owns\chi_k \mapsto r^{1/6} \, \chi_{\left[\frac{k-2r}{r^{1/3}},\frac{k+1-2r}{r^{1/3}}\right]}
\in L^2(\R)\,,
\quad k\in\Z\,,
\end{equation}
where $\chi_k\in\ell^2(\Z)$ is the characteristic function of the point
$k\in\Z$ and, similarly, the function on the right is the characteristic
function of a segment of length $r^{-1/3}$. Observe that this embedding is
isometric. Let $\J_r$ denote the kernel on $\R\times\R$ that is obtained from the
kernel $\J(\,\cdot\,,\,\cdot\,,r^2)$  on $\Z\times\Z$  using the embedding \eqref{e402}.  We shall establish
the following 

\begin{proposition}\label{p44} We have 
\begin{equation*}
[\J_r]_{[a,\infty)}\to[\A]_{[a,\infty)}\,, \quad r\to\infty \,,
\end{equation*}
in the trace norm for all $a\in\R$ uniformly on compact sets in $a$. 
\end{proposition}

This proposition immediately implies Theorem \ref{t4} as follows

\begin{proof}[Proof of Theorem \ref{t4}]
Consider the left-hand side of \eqref{e401} and choose for each  $a_i$ 
a pair of functions  $k_i^-(r),k_i^+(r)\in\Z$ such that
\begin{equation*}
\frac{k_i^-(r)-2r}{r^{1/3}}=a_i^-(r) \le a_i  \le  a_i^+(r) = \frac{k_i^+(r)-2r}{r^{1/3}} 
\end{equation*}
and $a_i^-(r),a_i^+(r)\to a_i$ as $r\to \infty$. Then, on the one hand, the 
probability in left-hand side of \eqref{e401} lies between the corresponding
probabilities for $a_i^-(r)$ and $a_i^+(r)$. On the other hand, the 
probabilities for $a_i^-(r)$ and $a_i^+(r)$ can be expressed in the form
\eqref{e41} for the kernel $\J_r$ and by Propositions \ref{p44} and 
continuity of the Airy kernel  they converge
to one and same limit given by the Airy kernel as $r\to\infty$.
\end{proof}

Now we get to the proofs of Propositions \ref{p41} and \ref{p44}  which will 
require some computations. Recall that 
the Airy function can be expressed in terms of  Bessel functions as follows 
\begin{equation}
A(x)=\cases \frac 1\pi \sqrt{\frac x3}K_{\frac 13}\left(\frac 23 x^{\frac 
32}\right)\,,&x\ge 0\,,\\
\frac{\sqrt{|x|}} 3\left[J_{\frac 13}\left(\frac 23 |x|^{\frac 
32}\right)+J_{-\frac 13}\left(\frac 23 |x|^{\frac 32}\right)\right]\,,&x\le 0\,,
\endcases
\end{equation}
see Section 6.4 in \cite{W}. 
Also recall that 
\begin{equation}\label{e42} 
A(x) \sim \frac1{2 x^{1/4} \sqrt{\pi}}\, e^{-\frac32 x^{3/2}}\,, \quad x\to+\infty \,,
\end{equation}
see, for example,  the formula 7.23 (1) in \cite{W}.

\begin{lemma}\label{l45} For any $x\in\R$ we have 
\begin{equation}\label{e43}
\left|r^{\frac 13}J_{2r+xr^{\frac 13}}(2r)- A(x)\right|= O(r^{-\frac 13})\,, \quad r\to\infty\,,
\end{equation}
moreover, the constant in  $O(r^{-\frac 13})$ is uniform in $x$ on compact subsets of $\R$. 
\end{lemma}

\begin{proof}
Assume first that $x\ge 0$. We denote 
\begin{equation*}
\nu=2r+xr^{\frac 13},\quad \alpha=\operatorname{arccosh}\left(1+xr^{-\frac 
23}/2\right)\ge 0.
\end{equation*}
It will be convenient to use the following notation
\begin{equation*}
P=\nu(\tanh\alpha-\alpha),\quad
Q=\frac \nu3 \tanh^3\alpha.
\end{equation*}
The formula 8.43(4) in \cite{W} reads
\begin{equation}\label{e44} 
J_{\nu}(2r)=\frac {\tanh\alpha}{\pi\sqrt{3}}e^{P+Q}K_{\frac 
13}\left(Q\right)+\frac{3\g_1}{\nu}\,e^{P}
\end{equation}
where $|\g_1|<1$. We have the following estimates as $r\to+\infty$
\begin{align*}
&\alpha=x^{\frac 12}r^{-\frac 13}+O(r^{-1}),\\
&\tanh\alpha=\alpha+O(\alpha^3)=x^{\frac 12}r^{-\frac 13}+O(r^{-1}),\\
&P+Q=\nu\cdot O(\alpha^5)=O(r^{-\frac 23}),\quad
e^{P+Q}=1+O(r^{-\frac 23}),\\
&Q=\frac 13\,\left(2r+xr^{\frac 13}\right)\left(x^{\frac 32}r^{- 1}+O(r^{-\frac 
43})\right)=\frac{2x^{\frac 32}}3+O(r^{-\frac 13}),\\
&K_{\frac 13}\left(Q\right)=K_{\frac 13}\left(\frac {2x^{\frac 
32}}{3}\right)+O(r^{-\frac 13}),\\
&P\le 0,\quad\frac{3\g_1}{\nu}\,e^{P}=O(r^{-1}).
\end{align*}
Substituting this into \eqref{e44}, we obtain the claim \eqref{e43}  for $x\ge 0$.

Assume now that $x\le 0$. Denote
\begin{equation*}
\nu=2r+xr^{\frac 13},\quad \beta=\operatorname{arccos}\left(1+xr^{-\frac 
23}/2\right)\ge 0,\quad y=|x|.
\end{equation*}
Introduce the notation
\begin{equation*}
\wt P=\nu(\tan\beta-\beta),\quad
\wt Q=\frac \nu 3\tan^3\beta\,.
\end{equation*}
The formula 8.43 (5) in \cite{W} reads
\begin{multline}\label{e45}
J_{\nu}(r)=\frac 13 \tan \beta\cos\left(\wt P-\wt Q\right)\left[J_{-\frac 
13}\left(\wt Q\right)+J_{\frac 13}\left(\wt Q\right)\right]\\
+\frac 1{\sqrt{3}} \tan \beta\sin\left(\wt P-\wt Q\right)\left[J_{-\frac 
13}\left(\wt Q\right)-J_{\frac 13}\left(\wt 
Q\right)\right]+\frac{24\g_2}{\nu}
\end{multline}
where $|\g_2|<1$. Again we have the estimates  as $r\to+\infty$
\begin{align*}
&\beta=y^{\frac 12}r^{-\frac 13}+O(r^{-1}),\\
&\tan\beta=\beta+O(\beta^3)=y^{\frac 12}r^{-\frac 13}+O(r^{-1}),\\
&\wt P-\wt Q=\nu\cdot O(\beta^5)=O(r^{-\frac 23}),\\
&\cos\left(\wt P-\wt Q\right)=1+O(r^{-\frac 43}),\quad
\sin\left(\wt P-\wt Q\right)=O(r^{-\frac 23}),\\
&\wt Q=\frac 13\left(2r-yr^{\frac 13}\right)\left(y^{\frac 32}r^{- 1}+O(r^{-\frac 
43})\right)=\frac{2y^{\frac 32}}3+O(r^{-\frac 13}),\\
&J_{\pm\frac 13}\left(\wt Q\right)=J_{\pm\frac 13}\left(\frac {2y^{\frac 
32}}{3}\right)+O(r^{-\frac 13}).
\end{align*}
These estimates after substituting into \eqref{e45} produce \eqref{e43}  for $x\le 0$. 
\end{proof}

\begin{lemma}\label{l46} There exist $C_1,C_2,C_3,\ep>0$ such that for any $A>0$ and $s>0$ we have
\begin{alignat}{2}
\left|J_{r+Ar^{\frac 13}+s}(r)\right|&\le
C_1\,{r^{-\frac13}}{\exp\left(-C_2\left(A^{\frac 32}+sA^{\frac
12}r^{-\frac 13}\right)\right)},\quad &s\le \varepsilon r \,, \label{e46a} \\
\left|J_{r+Ar^{\frac 13}+s}(r)\right|&\le
\exp\left(-C_3(r+s)\right), \quad &s\ge \varepsilon r\,, \label{e46b}
\end{alignat}
for all $r\gg 0$.
\end{lemma}

\begin{proof} 
First suppose  that $s\le\varepsilon r$. Set $\nu=r+Ar^{\frac 13}+s$. We shall use 
\eqref{e44}  with $\alpha=\operatorname{arccosh}(\nu/r)$. Provided $\varepsilon$ 
is chosen 
small enough and $r$ is sufficiently large,   $\alpha$ will be close to $0$ and we
will be able to use Taylor expansions.
For $r\gg 0$ we have
\begin{equation*}
\alpha=\operatorname{arccosh}(1+Ar^{-\frac 23}+sr^{-1})\ge
\const\,(Ar^{-\frac 23}+sr^{-1})^{\frac 12} \,,
\end{equation*}
and, similarly,
\begin{equation*}
-P=\nu(\alpha-\tanh\alpha)\ge \const\,(A+sr^{-\frac 13})^{\frac 32}\,.
\end{equation*}
Since the function $x^{\frac 32}$ is concave, we have 
$$ 
-P\ge \const\,(A^{\frac 32}+sA^{\frac 12}r^{-\frac 13})\,.
$$
The constant here is strictly positive.

Since $K_{\frac 13}(x)\le \const\, x^{-\frac 12} e^{-x}$ (see, for example, the 
formula 
7.23 (1) in \cite{W})  we obtain
\begin{multline*}
\tanh\alpha \,e^{P+Q}K_{\frac 13}\left(Q\right)\le \const\, 
\frac{e^{P}}{\sqrt{\nu\tanh\alpha}}\\ \le \frac{\const}  {r^{\frac 
13}}
\exp\left(-\const\,\left(A^{\frac 32}+sA^{\frac 12}r^{-\frac 13}\right)\right)\,,
\end{multline*}
where we used that  $\tanh\alpha\ge \const\, r^{-\frac 13}$.  Finally, we note that
\begin{equation*}
\frac {e^P}{\nu}\le \frac1r\, \exp\left({-\const\,\left(A^{\frac 32}+sA^{\frac 
12}r^{-\frac 13}\right)}\right),
\end{equation*}
and this completes the proof of \eqref{e46a}.

The estimate \eqref{e46b} follows directly from the formulas 
8.5 (9), (4), (5) in \cite{W}. 
\end{proof}

\begin{lemma}\label{l47} For any $\delta>0$ there exists such $M>0$ that for all 
$x,y>M$ and large enough $r$
\begin{equation*}
\left|\J\left(2r+xr^{\frac 13},2r+yr^{\frac 13}, {r^2}\right)\right|<\delta 
r^{-\frac 13}.
\end{equation*}
\end{lemma}

\begin{proof} From  \eqref{e214} we have 
\begin{equation}\label{e47}
\J\left(2r+xr^{\frac 13},2r+yr^{\frac 13}, {r^2} \right)=\sum_{s=1}^\infty 
J_{2r+xr^{\frac 13}+s}(2r)\,J_{2r+yr^{\frac 13}+s}(2r).
\end{equation}
Let us split the sum in \eqref{e47}  into two parts
\begin{equation*}
{\sum}_1=\sum_{l\le \varepsilon r}\,, \quad  
{\sum}_2=\sum_{l> \varepsilon r}\,,
\end{equation*}
that is, one sum for 
$l\le\varepsilon r$ and the other  for $l> \varepsilon r$, and apply
Lemma \ref{l46} to these two sums.  Note that $2r$ 
here corresponds to $r$ in Lemma \ref{l46};  this produces factors of 
$2^{\frac 13}$ and does not affect the estimate. 

Let the $c_i$'s stand for some positive constants not depending on 
$M$. From \eqref{e46a} we obtain the following estimate for the first sum 
$$
{\sum}_1 \le  c_1\, r^{-\frac 23}\,
\exp\left(-c_2\,M^{\frac 32}\right) \, \sum_{s=1}^{[\varepsilon r]} q^s 
$$
where
$$
q=\exp\left(-c_2M^{\frac 12}r^{-\frac 13}\right)\,, \quad 0<q<1\,.
$$
Therefore,
$$
{\sum}_1 \le \frac{c_1\, r^{-\frac 23}\,
\exp\left(-c_2\,M^{\frac 32}\right) }{1-q} \le 
r^{-\frac 13}\cdot c_3\exp(-c_2M^{\frac 32})M^{-\frac 12} \,.
$$
We can choose $M$ so that $c_3\exp(-c_2M^{\frac 32})M^{-\frac 12}<\delta/2$.

For  the second sum we use \eqref{e46b} and obtain  
\begin{equation*}
{\sum}_2 \le \sum_{s\ge \varepsilon r} \exp(-c_4 (r+s))\le  c_5\exp(-c_4r).
\end{equation*}
Clearly, this is less than $\delta r^{-\frac 13}/2$ for $r\gg0$.
\end{proof}

\begin{proof}[Proof of Proposition \ref{p41}] 
As shown in \cite{CL,TW}, the Airy kernel has the following integral 
representation
\begin{equation}\label{e48}
\A(x,y)=\int_{0}^\infty A(x+t)A(y+t)dt.
\end{equation}
The formula \eqref{e47} implies that
for any integer $N>0$
\begin{multline}\label{e49}
\J\left(2r+xr^{\frac 13},2r+yr^{\frac 13},{r^2} \right)=
\sum_{s=1}^N J_{2r+xr^{\frac 13}+s}(2r)\,J_{2r+yr^{\frac 13}+s}(2r)\\ 
+\J\left(2r+xr^{\frac 13}+N,2r+yr^{\frac 13}+N,{r^2} \right)\,.
\end{multline}
Let us fix $\delta>0$ and pick $M>0$ according to Lemma \ref{l47}.  Since, by assumption,
$x$ and $y$ lie in compact set of $\R$, we can fix $m$ such that
$x,y \ge m$. Set
$$
N=[(M-m+1)\, r^\frac13] \,.
$$
Then the inequalities 
\begin{equation*}
x+Nr^{-\frac 13}>M,\quad y+Nr^{-\frac 13}>M
\end{equation*}
are satisfied for all $x,y$ in our compact set and  Lemma \ref{l45} applies to
the sum in \eqref{e49}. We obtain 
\begin{multline*}
\left|r^{\frac 23}\sum_{s=1}^{N} J_{2r+xr^{\frac 13}+s}(2r)\,J_{2r+yr^{\frac 
13}+s}(2r)-\right.\\
\left. \sum_{s=1}^{N}
A(x+sr^{-\frac 13})A(y+sr^{-\frac 13})\right|=O(1)
\end{multline*}
because  the number of summands is $N=O(r^{\frac 13})$ and $A(x)$ is bounded on 
subsets of $\R$ which are bounded from below.
Note that 
\begin{equation*}
r^{-\frac 13} \sum_{s=1}^{N}
A(x+sr^{-\frac 13})A(x+sr^{-\frac 13})
\end{equation*}
is a Riemann integral  sum for the integral
\begin{equation*}
\int\limits_{0}^{M-m+1} A(x+t)A(y+t)\,dt,
\end{equation*}
and it converges to this integral as $r\to+\infty$. Since the absolute value of 
the second term in the right-hand side  of \eqref{e49}  does not exceed $\delta r^{-\frac 13}$ by 
the choice of $N$, we get
\begin{equation*}
\left|r^{\frac 13}\J\left(2r+xr^{\frac 13},2r+yr^{\frac 13},{r^2} 
\right)-\int\limits_{0}^{M-m+1} A(x+t)A(y+t)dt\right|\le \delta+o(1)
\end{equation*}
as $r\to+\infty$, and this estimate is uniform on compact sets. Now let  $\delta\to 
0$ and  $M\to+\infty$. By \eqref{e42} the
integral \eqref{e48} converges uniformly in $x$ and $y$  on compact sets and 
we obtain the claim of the proposition. 
\end{proof}

\begin{proof}[Proof of Proposition \ref{p44}]
It is clear that Proposition \ref{p41}  implies the convergence 
of $[\J_r]_a$ to $[\A]_a$ in 
the weak operator topology. Therefore, by Proposition \ref{p56}, it remains to prove that 
$\tr[\J_r]_a\to\tr[\A]_a$ as $r\to+\infty$. We have
\begin{equation*}
\tr[\J_r]_a=\sum_{k=[2r+ar^{\frac13}]}^\infty\J(k,k;{r^2})+o(1)\,,
\end{equation*}
where the $o(1)$ correction comes from the fact that $a$ may not be a number
of the form $\frac{k-2r}{r^{1/3}}$,  $k\in\Z$.  By \eqref{e47} we have 
\begin{equation}\label{e403}
\sum_{k=[2r+ar^{\frac13}]}^\infty\J(k,k;{r^2})
=\sum_{l=1}^\infty l\left(J_{[2r+ar^{\frac13}]+l}(2r)\right)^2 \,.
\end{equation}
Similarly, 
\begin{equation}\label{e404}
\tr[\A]_a=\int_a^{\infty}\A(s,s)ds=\int_0^\infty t(A(a+t))^2dt\,.
\end{equation}

Since we already established  the uniform convergence of kernels on compact sets, it is 
enough to show that the both \eqref{e403} and \eqref{e404} go to zero as $a\to+\infty$ and 
$r\to+\infty$. For the Airy kernel this is clear from \eqref{e42}. For the kernel $\J_r$
it is equivalent to the following statement: for any $\delta>0$ there 
exists $M_0>0$ such that for all $M>M_0$ and large enough $r$
we have 
\begin{equation}\label{e410}
\left|\sum_{l=1}^\infty l \, J_{2r+Mr^{\frac 13}+l}^2(2r)\right|<\delta \,. 
\end{equation}
We shall employ Lemma \ref{l46} for $A=M$. Again, we split the sum 
 in \eqref{e410} in two parts
\begin{equation*}
{\sum}_1=\sum_{l\le \varepsilon r}\,, \quad  
{\sum}_2=\sum_{l> \varepsilon r}\,.
\end{equation*}
For the first sum Lemma \ref{l46} gives
\begin{equation*}
{\sum}_1\le 
c_1r^{-\frac 23}\exp\left(-c_2M^{\frac 32}\right)\sum_{l\le [\varepsilon 
r]}l\, q^l\,,
\end{equation*}
where
\begin{equation*}
q=\exp\left(-c_2M^{\frac 12}r^{-\frac 13}\right)\,, \quad 0<q<1\,, 
\end{equation*}
and the $c_i$'s are some  positive constants that do not depend on $M$. 
Since $\sum l\, q^l = q (1-q)^{-2}$ we obtain
\begin{equation*}
{\sum}_1 \le 
c_1r^{-\frac 23}\exp\left(-c_2M^{\frac 32}\right) \, \frac{q}{(1-q)^2} \le
c_3\frac{\exp\left(-c_2M^{\frac 32}\right)}{M} \,. 
\end{equation*}
This can be made arbitrarily small by taking $M$ sufficiently large.

For the other part of the sum we have the estimate
\begin{equation*}
{\sum}_2 \le
\sum_{l> \varepsilon r}l \exp(-c_4(r+l))
\end{equation*}
which, evidently, goes to zero as $r\to+\infty$.
\end{proof}

\subsection{Depoissonization and proof of Theorem \ref{t3}}\label{s43}

Fix some $m=1,2,\dots$ and denote by $F_n$ the distribution
function of $\la_1,\dots,\la_m$ under the Plancherel measure $M_n$
\begin{equation*}
F_n(x_1,\dots,x_m)=M_n\left(\left\{\la\,\left|\,
\la_i< x_i\,, \  1\le i\le 
m\right.\right\}\right) \,.
\end{equation*}
Also, set
\begin{equation*}
F(\ix,x)=e^{-\ix} \sum_{k=0}^\infty\frac {\ix^k}{k!}\, F_k(x).
\end{equation*}
This is the distribution function corresponding to the measure $M^\ix$.

The measures $M_n$ can be obtained as distribution at time $n$ of
a certain random growth process of a Young diagram, see e.g.\ \cite{VK2}. This implies
that
\begin{equation*}
F_{n+1}(x)\le F_n(x)\,, \quad x\in\R^m \,.
\end{equation*}
Also, by construction, $F_n$ is monotone in $x$ and similarly 
\begin{equation}\label{e4101}
F(\ix,x) \le F(\ix,y)\,, \quad x_i \le y_i\,, \quad i=1,\dots,m \,. 
\end{equation}
We shall use these monotonicity properties 
together with the following lemma.  

\begin{lemma}[Johansson, \cite{J}]\label{l48}
 There exist constants $C>0$ and $n_0>0$ such that for 
any nonincreasing sequence $\{b_n\}_{n=0}^\infty\subset[0,1]$ 
\begin{equation*}
1\ge b_0\ge b_1\ge b_2\ge b_3\ge\dots\ge 0,
\end{equation*}
and its exponential generating function 
\begin{equation*}
B(\ix)=e^{-\ix} \sum_{k=0}^\infty\frac {\ix^k}{k!}\cdot b_k
\end{equation*}
we have for all $n>n_0$ the following inequalities:
\begin{equation*}
B(n+4\sqrt{n\ln n})-\frac C{n^2}\le b_n\le B(n-4\sqrt{n\ln n})+\frac C{n^2}.
\end{equation*}
\end{lemma}

This  lemma implies that for all $x\in \R^m$ 
\begin{equation}\label{e411}
F(n+4\sqrt{n\ln n},x)-\frac C{n^2}\le F_n(x)\le F(n-4\sqrt{n\ln n},x)+\frac 
C{n^2}\,.
\end{equation}
Set 
\begin{equation*}
\bo=(1,\dots,1)\,.
\end{equation*}
Theorem \ref{t4} asserts  that 
\begin{equation}\label{e412}
F\left(\ix, 2{\ix}^\frac12\, \bo + \ix^{\frac 16} \, x\right)\to F(x),\quad \ix\to +\infty,\quad x\in 
\R^m,
\end{equation}
where $F(x)$ is the corresponding distribution
function for the Airy ensemble. Note that $F(x)$ is continuous. 

 Denote $n_\pm=n\pm 4\sqrt{n\ln n}$. Then for $i=1,\dots,m$
\begin{equation*}
2n^\frac12_\pm +n^\frac16_\pm\, x_i =2n^\frac12  +n^{\frac 16}\, x_i +O((\ln n)^{1/2}) \,.
\end{equation*}
Hence, for any $\varepsilon>0$ and all sufficiently large $n$ we have 
\begin{equation*}
 2n^\frac12_+ +n_+^{\frac 16} \, (x_i-\varepsilon)
\le
2n^\frac12 +n^{\frac 16}\, x_i 
\le 2n^\frac12_- +n_-^{\frac 16} \, (x_i+\varepsilon) \,,
\end{equation*}
for $i=1,\dots,m$. By \eqref{e4101} this implies that 
\begin{align*}
F\left(n_+,2n^\frac12\,\bo+n^{\frac 16}\, x\right)&\ge 
F\left(n_+,2n^\frac12_+\, \bo+n_+^{\frac 16} \, \left(x-\varepsilon\, \bo\right)\right)\\
F\left(n_-,2n^\frac12\,\bo+n^{\frac 16}\, x\right)&\le 
F\left(n_-,2n^\frac12_-\, \bo+n_-^{\frac 16} \, \left(x+\varepsilon\, \bo\right)\right) \,. 
\end{align*} 
{}From this and \eqref{e411}  we obtain
\begin{multline*}
F\left(n_+,2n^\frac12_+\, \bo+n_+^{\frac 16} \, \left(x-\varepsilon\, \bo\right)\right)-\frac C{n^2}\\ \le 
F_n\left(2n^\frac12\,\bo+n^{\frac 16}\, x \right)\le \\  
F\left(n_-,2n^\frac12_-\, \bo+n_-^{\frac 16} \,\left(x+\varepsilon\, \bo\right)\right)+\frac C{n^2} \,. 
\end{multline*} 
{}From this and  \eqref{e412} we conclude that 
\begin{equation*}
F\left(x-\varepsilon\,\bo\right)+o(1)\le F_n\left(2n^\frac12\,\bo+n^{\frac 16}\, x \right)\le 
F\left(x+\varepsilon\, \bo\right)+o(1)
\end{equation*}
as $n\to\infty$. Since $\varepsilon>0$ is arbitrary and $F(x)$ is continuous 
we obtain
\begin{equation*}
F_n \left(2n^\frac12\,\bo+n^{\frac 16}\, x \right)\to F(x),\quad n\to \infty, \quad x\in \R^m,
\end{equation*}
which is the statement of Theorem \ref{t3}. 

\appendix 

\section{General properties of determinantal point processes}\label{s5}

In this Appendix, we collected some necessary facts about
determinantal point processes, their correlation functions, 
Fredholm determinants, 
and convergence of trace class operators. 

Let $\X$ be a countable set, let $\CX=2^{\X}$ be the set of subsets
of $\X$ and denote by $\CXf\subset\CX$ the set of finite subsets of $\X$.
We call elements of $\CX$ configurations. 
Let $L$ be a kernel on $\X$, that is, a function on $\X\times\X$
also viewed as a matrix of an operator in $H=\ell^2(\X)$. 

By a \emph{determinantal point process} on $\X$ (in \cite{DVJ}
such processes are called fermion point processes) 
 we mean a probability measure on
$\CXf$ such that
$$
\Prob(X)=\frac{\det\big[L(x_i,x_j)\big]_{x_i\in X}}{\det(1+L)}\,, \quad X\in\CXf\,.
$$
Here the determinant in the numerator is the usual determinant of linear
algebra, whereas
the determinant in the denominator is, in general, a Fredholm determinant.
Some sufficient conditions under which $\det(1+L)$ makes sense
are described in the following subsection.

\subsection{Fredholm determinants and determinantal processes}

Let $H$ be a complex Hilbert space, $\Lc(H)$ be the algebra of bounded
operators in $H$, and $\Lc_1(H)$, $\Lc_2(H)$ be the ideals of trace class
and Hilbert--Schmidt operators, respectively.

Assume we are given a splitting $H=H_+\oplus H_-$. 
According to this splitting, write operators $A\in\Lc(H)$ in block form 
$A=\bmatrix A_{++} & A_{+-}\\ A_{-+} & A_{--}\endbmatrix$, where 
\begin{align*}
&A_{++}: H_+\to H_+, \quad A_{+-}: H_-\to H_+, \\
&A_{-+}: H_+\to H_-, \quad A_{--}: H_-\to H_-\,.
\end{align*}
The algebra $\Lc(H)$ is equipped with a natural $\Z_2$-grading. Specifically,
given $A$, its even part $A_{\textup{even}}$ and odd part $A_{\textup{odd}}$
are defined as follows
\begin{equation*}
A_{\textup{even}}=\bmatrix A_{++} & 0\\ 0 & A_{--}\endbmatrix, \qquad
A_{\textup{odd}}=\bmatrix 0 & A_{+-}\\ A_{-+} & 0\endbmatrix \,.
\end{equation*}

Denote by $\LH$ the set of operators $A\in \Lc(H)$ such that $A_{\textup{even}}$ is in the
trace class $\Lc_1(H)$ while $A_{\textup{odd}}$ is in the Hilbert--Schmidt class 
$\Lc_2(H)$. It is readily seen that $\LH$ is an algebra. We endow it with the
topology induced by the trace norm on the even part and the Hilbert--Schmidt
norm on the odd part. 

It is well known that the determinant $\det(1+A)$ makes sense if $A\in\Lc_1(H)$. 
It can be characterized as the only function which is continuous in $A$ with respect
 to the trace norm $\Vert A\Vert_1=\tr\sqrt{AA^*}$ and which coincides with the usual 
determinant when $A$ is a finite--dimensional operator. See, e.g., \cite{Si}.

\begin{proposition}\label{p51} The function $A\mapsto \det(1+A)$ admits a unique extension to $\LH$,
 which is continuous in the topology of that algebra.
\end{proposition}

\begin{proof} For $A\in\LH$, set
\begin{equation}\label{e51}
\det(1+A)=\det((1+A)e^{-A})\cdot e^{\tr A_{\textup{even}}} 
\end{equation}
As is well known (e.g., \cite{Si}), 
\begin{equation*}
A\mapsto (1+A)e^{-A}-1
\end{equation*}
is a continuous map from $\Lc_2(H)$ to $\Lc_1(H)$.  Next, $A\mapsto \tr A_{\textup{even}}$
evidently is a continuous function on $\LH$. Consequently, \eqref{e51}  is well 
defined and is a continuous function on $\LH$. 
When $A\in\Lc_1(H)$, \eqref{e51} agrees with the conventional definition, because then 
\begin{equation*}
\det((1+A)e^{-A})\cdot e^{\tr A_{\textup{even}}}=\det(1+A)e^{-\tr A+\tr A_{\textup{even}}}=
\det(1+A).
\end{equation*}
This concludes the proof. \end{proof}

\begin{corollary}\label{c52} If $\{P_n\}$ is an ascending sequence of projection operators in $H$ such that $P_n\to1$ strongly, then 
\begin{equation*}
\det(1+A)=\lim_{n\to\infty}\det(1+P_n A P_n).
\end{equation*}
\end{corollary} 

\begin{proof} Indeed, $P_nAP_n$ approximates $A$ in the topology of $\LH$. \end{proof}

\begin{corollary}\label{c53} If $A,B\in\LH$ then
\begin{equation*}
\det(1+A)\, \det(1+B)=\det((1+A)(1+B)).
\end{equation*}
\end{corollary}

\begin{proof} Indeed, this is true for finite--dimensional $A,B$, and then we use 
the continuity argument. 
\end{proof}

In our particular case, the splitting of $H=\ell^2(\X)$
will come from a splitting of $\X=\X_+\sqcup\X_-$ into
two complementary subsets as follows
$$
H_\pm=\ell^2(\X_\pm)\,.
$$
An operator $L$ in $H$ will be viewed as an infinite matrix
whose rows and columns are indexed by elements of $\X$. Given $X\subset\X$, we
denote by $L_X$ the corresponding finite submatrix in $L$. 

\begin{proposition}\label{p54} If $L\in\LH$ then
\begin{equation}\label{e52}
\sum_X\det L_X=\det(1+L),
\end{equation}
where summation is taken over all finite subsets $X\subset\X$ including the
empty set with understanding that $\det L_\varnothing=1$. 
\end{proposition} 

\begin{proof} Given a finite subset $Y\subset \X$, we assign to it, in the
 natural way, a projection operator $P_Y$. Then, by elementary linear algebra,
\begin{equation*}
\sum_{X\subseteq Y}\det L_X=\det(1+P_YLP_Y).
\end{equation*}
Assume $Y$ becomes larger and larger, so that in the limit it 
covers the whole $\X$. Then the left-hand side tends to the left-hand side
 of \eqref{e52}. 
On the other hand the right-hand side tends to $\det(1+L)$ by Corollary \ref{c52}.  
\end{proof}

\begin{remark}
Suppose that 
 $L=\bmatrix 0 & A\\ -A^* & 0\endbmatrix$, where $A$ is of Hilbert--Schmidt
class. Then $L\in\LH$. It is readily seen that $\det L_X\ge0$ for all $X$, and it is  worth noting that $\det L_X=0$ unless $|X_+|=|X_-|$. 
By Proposition \ref{p54}, we can define a probability measure 
on finite subsets $X$ of $\X$ by
\begin{equation*}
\Prob(X)=\frac{\det L_X}{\det(1+L)}\,, \quad X\in\CXf\,.
\end{equation*}
\end{remark} 

\subsection{Correlation functions of determinantal processes}

Given $X\in\CXf$, let $\rho(X)$ be the probability that a 
random configuration  contains $X$, that is,
$$
\rho(X)=\Prob\left(\left\{Y\in\CXf, X\subset Y\right\}\right)\,.
$$
We call $\rho(X)$ the \emph{correlation functions}. The fundamental
fact about determinantal point processes is that their
correlation functions have again a determinantal form. 

\begin{proposition}\label{p55} Let $L$ be as above and set $K=L(1+L)^{-1}$. 
Then $\rho(X)=\det K_X$. 
\end{proposition}

\begin{proof} We follow the argument in \cite{DVJ}, Exercise 5.4.7.  Let $f(x)$ be an arbitrary function on $\X$ such that $f(x)=1$ for all but a finite number of $x$'s. Form the probability generating functional:
\begin{equation*}
\Phi(f)=\sum_X \prod_{x\in X}f(x)\cdot \Prob(X).
\end{equation*}
Then, viewing $f$ as a diagonal matrix, we get
\begin{equation*}
\Phi(f)=\frac{\sum_X\det(fL)}{\det(1+L)}
=\frac{\det(1+fL)}{\det(1+L)}\,, 
\end{equation*}
where the last equality is justified by Proposition \ref{p54} applied to the operator $fL$.

Now, set $g(x)=f(x)-1$, so that $g(x)=0$ for all but finitely many $x$'s. Then 
we can rewrite this relation as follows
\begin{equation*}
\Phi(f)= \frac{\det(1+fL)}{\det(1+L)}=\frac{\det(1+L+gL)}{\det(1+L)}
=\det(1+gK),
\end{equation*}
where the last equality follows by Corollary \ref{c53}.

Next, as $gK$ is in $\LH$ (it is even finite--dimensional), this can be rewritten as
\begin{equation*}
\Phi(f)=\sum_X\det((gK)_X)=\sum_X\prod_{x\in X}g(x)\cdot\det K_X.
\end{equation*}
On the other hand, by the very definition of $\Phi(f)$,
\begin{equation*}
\Phi(f)=\sum_X \prod_{x\in X}g(x)\cdot\rho(X).
\end{equation*}
This implies $\rho(X)=\det K_X$, as desired. 
\end{proof}

\begin{remark} If $L=\bmatrix 0 & A\\ -A^* & 0\endbmatrix$ then 
\begin{equation*}
K=\bmatrix AA^*(1+AA^*)^{-1} & (1+AA^*)^{-1}A\\
-(1+A^*A)^{-1}A^* & A^*A(1+A^*A)^{-1}\endbmatrix\,.
\end{equation*}
\end{remark}

\subsection{Complementation principle}\label{sA3}

In this section we discuss a simple but useful observation
which was communicated to us by S.~Kerov. 
Consider an arbitrary probability measure on $\CX$ such that
its correlation functions
$$
\rho(X)=\Prob\left(\left\{Y\in\CX, X\subset Y\right\}\right)\,, \quad X\in\CXf\,,
$$
have a determinantal form
$$
\rho(X)=\det\Big[K(x_i,x_j)\Big]_{x_i\in X}
$$
for some kernel $K$. 

Let $Z\subset \X$ be an arbitrary subset of $\X$. Consider the 
symmetric difference mapping
\begin{equation*}
\trZ: \CX\to \CX\,,\quad
Y\mapsto Y\tri Z \,,
\end{equation*}
which is an involution in $\CX$. Let $\Prob^\tri=\left(\trZ\right)_* \Prob$
be the image of our probability measure under $\trZ$
and let $\rho^\tri(X)$ be the correlation functions of the 
measure $\Prob^\tri$. Define a new kernel $K^\tri$ as follows. Let $Z'=\X\setminus Z$
be the complement of $Z$ and write the matrix $K$ in the block form
with respect to the decomposition $\X=Z'\sqcup Z$
$$
K_{Z'\sqcup Z}=
\begin{bmatrix} A & B \\ C & D 
\end{bmatrix}\,.
$$
By definition, set
$$
K^\tri_{Z'\sqcup Z}=
\begin{bmatrix} A & B \\ -C & 1-D 
\end{bmatrix}\,.
$$
We have the following
\begin{proposition}\label{pAc}
 $\rho^\tri(X)=\det\Big[K^\tri(x_i,x_j)\Big]_{x_i\in X}$\,.
\end{proposition}
\begin{proof}
Set $X_1=X\setminus Z$, $X_2=Z\setminus X$.
By the inclusion-exclusion principle we have
\begin{align*}
\rho^\tri(X)&=\Prob\left(\left\{Y\in\CX, X_1\subset Y, X_2\cap Y=\emptyset
\right\}\right) \\
&=\sum_{S\subset X_2} (-1)^{|S|} \rho(X_1 \cup S)\,.
\end{align*}
This alternating sum is easily seen to be identical to
the expansion of $\det\Big[K^\tri(x_i,x_j)\Big]_{x_i\in X}$ by
linearity using 
$$
\begin{bmatrix} A & B \\ -C & 1-D 
\end{bmatrix} = 
\begin{bmatrix} A & B \\ -C & -D 
\end{bmatrix} + \begin{bmatrix} 0 & 0 \\ 0 & 1 
\end{bmatrix}
$$
\end{proof}

\subsection{Convergence of trace class operators}

Let $K_1,K_2,\dots$ and $K$ be Hermitian nonnegative operators in
$\Lc_1(H)$. The following proposition is a special case of Theorem 2.20 in the book
\cite{Si2} (we are grateful to P.~Deift for this reference). For the reader's convenience
we give a proof here. 

\begin{proposition}\label{p56} The following conditions are equivalent:

{\rm(i)} $\Vert K_n-K\Vert_1\to0$;

{\rm(ii)} $\tr K_n\to \tr K$ and $K_n\to K$ in the weak operator topology. 
\end{proposition}

First, prove a lemma:

\begin{lemma}\label{l57} Let $X=\bmatrix A & B\\ B^* & D \endbmatrix$ be a
nonnegative operator $2\times 2$ matrix. Then 
$\Vert B\Vert_1\le \sqrt{\tr A\cdot\tr D}$.
\end{lemma} 

\begin{proof}[Proof of Lemma] Without loss of generality one can assume that
the block $B$ is a nonnegative diagonal matrix,
$B=\textup{diag}(b_1,b_2,\dots)$. Write the blocks $A$ and $D$ as matrices,
too, and let $a_i$ and $d_i$ be their diagonal entries. Since
$X\ge0$, we have  $b_i^2\le a_id_i$ and therefore 
\begin{equation*}
\Vert B\Vert_1=\sum b_i\le\sum\sqrt{a_i d_i}\le 
\sqrt {\sum a_i\cdot\sum d_i}\le
\sqrt{\tr A\cdot\tr D}.
\end{equation*}
\end{proof}

\begin{proof}[Proof of Proposition \ref{p56}] Clearly, (i) implies (ii). To check the
converse claim, write $K$ in the block form, 
$K=\bmatrix A & B\\ B^* & D \endbmatrix$, where $A$ is of finite size
and $\tr D$ is small. Write all the $K_n$'s in the block
form with respect to the same decomposition of the Hilbert space, 
$K_n=\bmatrix A_n & B_n\\ B_n^* & D_n \endbmatrix$. Since $K_n\to K$
weakly, we have convergence of finite blocks, $A_n\to A$, which implies
$\tr A_n\to \tr A$. Since $\tr K_n\to \tr K$, we get $\tr D_n\to\tr
D$, so that all the traces $\tr D_n$ are small together with $\tr D$
provided that $n$ is large enough.

Write $K'=\bmatrix A & 0\\ 0 & 0 \endbmatrix$ and similarly for
$K_n$. Then 
\begin{equation*}
\Vert K_n-K\Vert_1\le 
\Vert K_n-K'_n\Vert_1+
\Vert K'_n-K'\Vert_1+
\Vert K'-K\Vert_1.
\end{equation*}
In the right-hand side, the first and the third summands are small because of the
lemma, while the second summand is small because it is equal to 
$\Vert A_n-A\Vert_1$. 
\end{proof}

\begin{proposition}\label{p58} The map $(A_1,\dots,A_n)\mapsto \det (I+\lambda_1
A_1+\dots+\lambda_n A_n)$ defines a continuous map from $(\Lc_1(H))^n$
to the algebra of entire functions in $n$ variables with the topology of
uniform convergence on compact sets.
\end{proposition} 
\begin{proof} The fact that $\det (I+\lambda_1 A_1+\dots+\lambda_n A_n)$ is
holomorphic in $\{\lambda_i\}$ for any trace class operators
$A_1,\dots,A_n$ is proved in \cite{Si}. The continuity of the map follows
from the inequality
\begin{equation*}
|\det(I+B)-\det(I+C)|\le\Vert B-C\Vert_1\exp(\Vert B\Vert_1+ \Vert
C\Vert_1+1)
\end{equation*}
which holds for any $B,C\in \Lc_1(H)$, see \cite{SS,Si2}.
\end{proof}

\end{document}